\documentclass[a4paper,11pt]{article}
\usepackage[latin1]{inputenc}
\usepackage[T1]{fontenc}
\usepackage{a4wide}
\usepackage{fontenc}
\usepackage{color}
\usepackage{amsmath}
\usepackage{graphicx}
\usepackage{array}
\usepackage{amsfonts}
\usepackage{amssymb}
\usepackage{euscript}
\usepackage{delarray}
\usepackage{latexsym}
\usepackage{enumerate}
\usepackage{float}
\usepackage{algorithmicx}
\usepackage{algorithm}
\usepackage{algpseudocode}
\usepackage{url}
\usepackage{breakcites}
\usepackage{mathtools}
\usepackage{dsfont}

\usepackage{caption}
\usepackage{subcaption}

\usepackage{multirow}
\usepackage{multicol}

\title{ANOVA decomposition of conditional Gaussian processes for sensitivity analysis with dependent inputs}
\author{Gaelle Chastaing \\ Loic Le Gratiet}

\newcommand{\etal}{\textit{et al.}}

\newcommand{\N}[2]{\mathcal{N}\left(#1, #2 \right)}
\newcommand{\XX}{\mathbf{X}}
\newcommand{\xx}{{\mathbf{x}}}
\newcommand{\tx}{{\mathbf{t}}}

\newcommand{\un}{\mathbf{1}_n}
\newcommand{\txx}{\tilde{\mathbf{x}}}
\newcommand{\R}{\mathbb{R}}
\newcommand{\PX}{P_\mathbf{X}}
\newcommand{\cov}{\mathrm{Cov}}
\newcommand{\EZ}{\mathbb{E}_Z}
\newcommand{\E}{\mathbb{E}}
\newcommand{\covZ}{\mathrm{Cov}_Z}
\newcommand{\D}{\mathbf{D}}
\newcommand{\TT}{\mathbf{T}}
\newcommand{\U}{\mathbf{U}}

\newcommand{\kk}{\mathbf{k}}
\newcommand{\K}{\mathbf{K}}
\newcommand{\MM}{\mathbf{M}}
\newcommand{\LL}{\mathbf{L}}
\newcommand{\Kn}{\mathbf{K}_n}
\newcommand{\kn}{\mathbf{k}_n}

\newcommand{\y}{\mathbf{z}_n}
\newcommand{\T}{{}^t}
\newcommand{\pX}{p_\mathbf{X}}
\newcommand{\pXu}{p_{\mathbf{X}_u}}

\newcommand{\Zn}{\boldsymbol{Z}_n}

\newcommand{\mmu}{\boldsymbol{\mu}}
\newcommand{\ttheta}{\boldsymbol{\theta}}
\newcommand{\pphi}{\boldsymbol{\phi}}

\newcommand{\GGamma}{\boldsymbol{\Gamma}}
\newcommand{\eepsilon}{\boldsymbol{\varepsilon}}

\newcommand{\ww}{\mathbf{w}}

\newcommand{\GP}{\mathrm{GP}}

\newtheorem{prop}{Proposition}

\begin{document}
\maketitle

\section*{Abstract}

Complex computer codes are widely used in   science   to model  physical systems. Sensitivity analysis aims to measure the contributions of the inputs  on the code  output variability.  An efficient tool to perform such analysis are the variance-based methods which have  been recently  investigated in the framework of dependent  inputs.  
One of their issue is that they  require  a large number of runs for the complex simulators. 
To handle it, a Gaussian process regression model may be used to approximate the complex code. 
In this work, we propose to decompose a Gaussian process into  a high dimensional  representation. This leads to  the definition of a  variance-based sensitivity measure well tailored for non-independent inputs.
We give a methodology to estimate these indices and to quantify their uncertainty.
Finally, the  approach  is illustrated on toy functions and on   a river flood model.\\

\textbf{Keywords:} Sensitivity analysis, dependent inputs, Gaussian process regression, functional decomposition, complex computer codes.

\section{Introduction}\label{introduction}

Many physical phenomena are investigated by complex models implemented in computer codes. Often considered as a black box function, a computer code calculates one or several output values which depend on input parameters. However, the code may depend on a very large number of incomes, that can be correlated among them.   Moreover, input parameters are also subject to many sources of uncertainty, attributed to errors of measurements or a lack of information.
These major flaws undermine the confidence a user have in the model. Indeed, the prediction given by the model may suffer from a large variability, leading to wrong conclusions. \\
To tackle these issues, the sensitivity analysis offers a series of methods and strategies that has been widely studied over the past decades~\cite{saltelli,saltelliglobal,cacuci}. Among the wide range of proposed methods, one could cite the class of global sensitivity analysis. Based on the assumption that the parameters implied in the model are randomly distributed, the global sensitivity analysis aims to identify and to rank the most contributive inputs to the response variability. 
One of the most popular global measure, the Sobol index, is based on a variance decomposition. Advanced by Hoeffding~\cite{hoeffding}, the model function can be uniquely decomposed as a sum of mutually orthogonal functional components when input variables are independent. Following this idea, Sobol constructs sensitivity measures by expanding the global variance into partial variances. Then, the Sobol index apportions the individual contribution of a set of inputs by the ratio between the partial variance depending on this set and the global variance~\cite{sobol}. \\
However, the construction of such measure relies on the assumption that input variables are independent. When incomes are dependent, the use of the Sobol index is not excluded, but it may lead to a wrong interpretation. Indeed, as underlined by Mara \etal~\cite{mara}, if inputs are not independent, the amount of the response variance due to a given factor may be influenced by its dependence to other inputs. In other word, as the Sobol index only depends on terms of variance, we ignore how it differentiates the inputs dependence from their interactions. From this perspective, the construction of a sensitivity measure that quantifies the uncertainty brought by dependent inputs becomes clear. \\
A solution is to use a  functional  decomposition to build a variance-based sensitivity index.
First, Xu \etal~\cite{xu} propose to decompose the partial variance of an input into a correlated and an uncorrelated contribution under the hypothesis that the effect of each parameter on the response is linear. The authors learn these contributions by successive linear regressions. To improve this approach, Li \etal~\cite{li} propose to approximate the model function by a High Dimension Model Representation (HDMR), that consists of a sum of functional components of low dimensions~\cite{hdmr}. They suggest to reconstruct each term via the usual basis functions (polynomials, splines,\dots). Then, they deduce the decomposition of the response variance as a sum of partial variances and covariances. Recently, Caniou \etal~\cite{caniouth}  suggest to build a HDMR by substituting the model function to a truncated polynomial chaos~\cite{wiener38} orthogonal with respect to the product of inputs marginal distribution. This choice is motivated by the fact that, when inputs are independent, the functional decomposition recovers the Hoeffding one, where each (unique) summand is expanded in terms of polynomial chaos~\cite{sudret}.\\
In a recent paper, Chastaing \etal~\cite{chastaing} revisit the Hoeffding decomposition in a different way. To tackle the problem of uniqueness of the components of the decomposition proposed by the previous approaches, the authors give a unique decomposition of the theoretical model. The main strength of the approach is that it is not based on surrogate modeling. Initiated by the pioneering work of Stone~\cite{stone}, they show that any regular function can be uniquely decomposed as a sum of hierarchically orthogonal component functions. This means that two of these summands are orthogonal whenever all variables included in one of the component are also involved in the other. The decomposition leads to the definition of a generalized sensitivity index involving variance and covariance components. Further, the same authors propose a numerical method of estimation~\cite{chastaing2}. \\
However, all these approaches suffer from two major flaws for time-consuming computer codes. First, the estimation of these measures is done by a regression method, which requires a very large number of model evaluations to be robust. Secondly, the number of decomposition components exponentially grows with the model dimension. In practice, we assume that only the low-order interaction terms contain the major part of the model behaviour. However, very few theoretical arguments confirm this assumption, and the truncation leads to an error of approximation that can be hardly controlled. 

To overcome the first issue, we surrogate the computer code with a Gaussian process regression model.  It is a non-parametric approach which considers that our prior knowledge about the code  can be modeled by a Gaussian process (GP)~\cite{santner,rasmussen}. These models are widely used in computer experiments  to surrogate a complex computer code from few of its outputs (\cite{SACKS89}). Further, the use of a GP model for sensitivity analysis is motivated by arguments given in the literature. For independent inputs, a natural approach is to substitute the model function by the posterior mean of a given GP~\cite{chen}. As this approach does not consider the posterior variance of the GP and thus the uncertainty of the surrogate modeling, 
 Oakley \& O'Hagan~\cite{oakley} substitute the initial model to a GP in the Sobol index. Then, the sensitivity index is given by the posterior mean of the Sobol index whereas the posterior variance measures its uncertainty. These two approaches are investigated and numerically compared in Marrel \etal~\cite{marrel}.\\
To handle with the second issue relative to the decomposition  truncation, we propose here to extend the work done by Durrande \etal~\cite{durrande} to the case of models with dependent inputs. In particular, we deal with GP specified by a covariance kernel that belongs to a special class of ANOVA kernels studied in~\cite{berlinet,durrande}. But instead of considering the posterior mean, we propose a functional decomposition of a GP distributed with respect to   the posterior distribution. Similarly to the work of Caniou \etal~\cite{caniouth}, the considered GP   is decomposed as a sum of processes indexed by increasing dimension input variables. This expansion is such that the summands are mutually orthogonal with respect to the product of the inputs marginal distributions. Thus, as we have accessed to the GP, and as we can deal with each term of its decomposition, we can easily deduce the sum of every other terms, so that a truncation error can not be produced.

Consequently, the GP development leads to the construction of sensitivity measures based on the decomposition of the global variance as a sum of partial variances and covariances. The difference with the use of the polynomial chaos is that GP are not intrinsically linked to the distribution of the input variables, as it is the case for polynomial chaos~\cite{cameron}.  Also, it should be noticed that our measure is a distribution which takes into account the uncertainty of the surrogate modeling. Further, we propose a numerical method to estimate our new defined measures. 
The procedure is experimented on several numerical examples. Furthermore, we study the asymptotic properties of the estimated measures.\\

The paper is organized as follows. In Section \ref{GPreg},  we introduce the first definitions and the main features of a GP. We also study a special case of covariance kernels. They will be used all along the article as the referenced kernels because they have good properties for sensitivity analysis. In Section \ref{anovadec}, we first remind the ANOVA decomposition proposed by Durrande \etal~\cite{durrande}. This expansion is done on the posterior mean of a GP. This introduces the decomposition of the conditional Gaussian processes we develop in this paper. After giving the advantages of such approach, we define in Section \ref{si} a sensitivity index well suited for models with dependent inputs. Furthermore, we describe a numerical procedure based on the Monte Carlo estimate to compute our new sensitivity index. At the end of Section \ref{si}, we study the convergence properties of the measure. Section \ref{applications} is devoted to numerical applications. The goal is to show the relevance of such a sensitivity index through several test cases. Further, we apply our methodology on a real-world problem.

\section{Gaussian process regression for sensitivity analysis}\label{GPreg}

In this section, we introduce the notation that will be used along the document. We  remind the basic settings on Gaussian process regression. The purpose is to build a fast approximation --- also called meta-model --- of  the input/output relation of the objective function. Then, we present a particular Gaussian process regression which is relevant to perform sensitivity analysis with dependent inputs.

\subsection{First definitions}

Let $(\Omega_\XX,\mathcal A_\XX,\mathbb P_\XX)$ be a probability space. Let $f$ be a measurable function of a random vector $\XX=(X_1,\cdots,X_p) \in \mathbb R^p$, $p \geq 1$, and defined as,

\[
\begin {array}{cccccc}
&(\Omega_\XX,\mathcal A_\XX,\mathbb P_\XX) & \rightarrow & (\R^p, \mathcal B(\R^p), \PX) & \rightarrow & (\R,\mathcal B(\R))\\
f: &\omega& \mapsto & \XX(\omega) & \mapsto & f(\XX(\omega)),
\end{array}
\]
where the joint distribution of $\XX$ is denoted by $\PX$. Further, we assume that $\PX$ is absolutely continuous with respect to the Lebesgue measure, and that $\XX$ admits a density $\pX$ with respect to the Lebesgue measure, i.e. $\pX d \xx = d \PX$.\\
Also, we assume that $f \in L_\R ^2(\R^p, \mathcal B(\R^p), \PX)$. We define the expectation with respect to $\PX$ as follows,

\[
\E(h(\XX))=\int_{\R^p} h(\xx) \pX(\xx) d\xx, \quad h\in L_\R ^2(\R^p, \mathcal B(\R^p), \PX).
\]
Further, $V(\cdot)=\E(\cdot-\E(\cdot))^2$ denotes the variance, and $\cov(\cdot,*)=\E[(\cdot-\E(\cdot))(*-\E(*))]$ the covariance with respect to the inputs distribution $\PX$.\\
The collection of all subsets of $\{1,\dots,p\}\setminus\{\emptyset\}$ is denoted by $S$. For $u \in S$ with $u=(u_1,\cdots,u_t)$, $|u|=t\geq 1$, the random subvector $\XX_u$ of $\XX$ is defined as $\XX_u:=(X_{u_1},\cdots,X_{u_t})$. 
The marginal density of $\XX_u$  is denoted by $\pXu$.

\subsection{Introduction to Gaussian process regression}\label{GPRequations}

For $\xx=(x_1,\cdots,x_p) \in \R^p$, we consider that the prior knowledge about $f(\xx)$ can be modeled by a zero-mean Gaussian Process (GP) $Z(\xx)$ defined on a probability space $(\Omega_Z, \mathcal{A}_Z,\mathbb P_Z)$ plus a known mean $m(\xx)$,
\begin{displaymath}
f(\xx) = m(\xx) + Z(\xx).
\end{displaymath}
From now, we denote by $\EZ$, $V_Z$ and $\covZ$ the expectation, variance and covariance with respect to $\mathbb P_Z$. 
A GP is completely specify by its mean $\EZ[Z(\xx)] $ and its covariance kernel  $k(\xx, \txx) = \covZ(Z(\xx), Z(\txx))$. Here, we consider a zero-mean GP, that can be written as 
\begin{equation}\label{gp1}
Z(\xx) \sim \mathrm{GP}(0, k(\xx, \txx)).
\end{equation}
Further, we denote by $ \D = \{\xx^1,\dots, \xx^n\}$, with $\xx^j \in \R^p$ for $j=1,\cdots,n$, the $n$-sample of observed inputs. We write the vector of centered outputs $\y := {}^t (f(\xx^1) - m(\xx^1) , \cdots,  f(\xx^n) - m(\xx^n) )$. Further, we consider the Gaussian random  vector $\Zn =  (  Z(\xx^1), \cdots ,  Z(\xx^n))$.  Notice that $(\xx^j)_{j=1,\dots,n}$ are generally not  sampled from the distribution $P_\XX$. Indeed, they usually come from a space-filling design procedure \cite{FAN06} in order to obtain good prediction accuracy.  

In the kriging theory, the aim is to use  the known values $\y  $   of $\Zn$ at points in $\D$ to predict $f(\xx)$. To perform such  prediction, we consider the conditional distribution $[f(\xx) | \Zn =\y]$. 
Standard results about Gaussian distribution gives that  this  conditional distribution   is given by 
\begin{equation}\label{predictivedistribution}
[f(\xx) | \Zn = \y ] = \mathrm{GP}\left(\mu(\xx), s^2(\xx, \txx) \right),
\end{equation}
with 
\begin{equation}\label{predictivemean}
\mu(\xx) = m(\xx) + \T \kn(\xx) \Kn^{-1}(\y - \mathbf{m}_n), 
\end{equation}
and 
\begin{equation}\label{predictivecovariance}
 s^2(\xx,\txx) = k(\xx, \txx) -  \T \kn(\xx) \Kn^{-1}\kn(\txx), 
\end{equation}
where $\kn(\xx) =  [k(\xx , \xx^j)]_{j=1,\dots,n}$, 
$ \mathbf{m}_n = {}^t (m(\xx^1) \cdots m(\xx^n))$ and~$ \Kn = [k(\xx^j, \xx^l)]_{j,l=1,\dots,n}$.\\

The mean $\mu(\xx)$ of the predictive distribution $[f(\xx) | \Zn =\y]$ is considered as the meta-model  for $f(\xx)$ and $s^2(\xx, \xx)$ represents   its  mean  squared error.
An important property of Gaussian process regression  is that the mean $\mu(\xx)$ interpolates the observations $\y$ and the variance $ s^2(\xx,\xx)$ equals zero at points in $\D$.

\subsection{Covariance kernel  for sensitivity analysis}

Certainly one of  the most important points of  Gaussian process regression is the choice of the covariance kernel $k(\xx, \txx) $,   for $\xx = (x_1, \dots, x_p)$, $\txx = (\tilde x_1, \dots, \tilde x_p)\in \mathbb R^p$,  of the unconditioned Gaussian process $Z(\xx) $ modeling the residual $f(\xx) - m(\xx)$.  We note that $k(\xx, \txx) $ must be   positive definite and we consider here  that $\sup_{\xx, \txx \in \R^p} k(\xx, \txx)  < \infty$. We choose in this paper a relevant class of kernels  for performing sensitivity analysis. They are built from Proposition \ref{prop1} \cite{durrande}.

\begin{prop}\label{prop1}
Let us consider a covariance kernel $\tilde{k}(\xx, \txx)$, $\xx, \txx \in \R^p$, such that $\tilde{k}_\xx: \xx \mapsto \tilde{k}(\xx, \txx)$ is in $L_1(\R^p)$  for all $\xx \in \R^p$ and  $\tilde{k}: (\xx, \txx)  \mapsto \tilde{k}(\xx, \txx)$ is in $L_1(\R^p \times \R^p)$. Then the following kernel $k(\xx, \txx) $  is a covariance kernel:
\begin{equation}\label{kernel0}
k(\xx, \txx) = \tilde{k}(\xx, \txx) - \frac{\int {\tilde{k}(\xx, \ww) \pX(d\ww) } \int{ \tilde{k}(\tilde \ww, \txx) \pX(d\tilde \ww)} }{\int\!\!\!  \int \tilde{k}(\ww, \tilde \ww) \pX(d\ww)  \pX(d\tilde \ww)  }.
\end{equation}
Furthermore, if we consider a  Gaussian process $Z(\xx) \sim  \mathrm{GP}\left(0, k(\xx, \txx) \right)$, then we have the following equality almost surely,
\begin{displaymath}
\int {Z( \xx) \pX(d\xx)} = 0.
\end{displaymath}
\end{prop}

%
%
From now and until the end of the article, we are interested by the following covariance kernel:
\begin{equation}\label{anovakernel}
k(\xx, \txx) = \sigma^2 \prod_{i=1}^p \left( 1 +  k_0^i(x_i, \tilde{x}_i)  \right), 
\end{equation}
where, following Proposition \ref{prop1}, for all $i=1,\dots,p$, we set:
\begin{equation}\label{k0i}
k_0^i(x_i, \tilde{x}_i) = \tilde{k}^i(x_i, \tilde{x}_i) - \frac{\int {\tilde{k}^i(x_i, w) p_{X_i}(w)dw } \int{ \tilde{k}^i(v, \tilde{x}_i) p_{X_i}(\tilde w)d\tilde w} }{\int \!\!\! \int \tilde{k}^i(w, \tilde w) p_{X_i}(w) p_{X_i}(\tilde w) dw  ~ d\tilde w  }, 
\end{equation}
where  $(\tilde{k}^i(x_i, \tilde{x}_i) )_{i=1,\dots,p}$ are given covariance kernels such that $\tilde{k}^i_w: \tilde w \mapsto \tilde{k}^i(w, \tilde w)$ is in $L_1(\R)$  for all $w \in \R$ and  $\tilde{k}^i: (w, \tilde w)  \mapsto \tilde{k}(w, \tilde w)$ is in $L_1(\R \times \R)$. A nice property of $k_{0 }^i(x_i, \tilde{x}_i)$,  $i=1,\dots,p$, is that it  is centered with respect to the marginal probability density function $p_{X_i}$. This feature is going to be exploited in Section \ref{anovadec}. \\
The choice of $\tilde{k}^i(x, \tilde x)$,  ${i=1,\dots,p}$,  is  of importance since it controls the regularity  in the $i^{\mathrm{th}}$ direction of the Gaussian process $Z(\xx)$  and thus the smoothness of the meta-model (see \cite{S99} and \cite{rasmussen}). For instance, for $m \in \mathbb N$, the partial derivative $\partial^m Z(\xx) / \partial^m x_i$ exists in mean square sense if and only if  the $2m$ derivative   of $\tilde{k}^i_0(x_i, \tilde{x}_i)$ exists at point $x_i = \tilde{x}_i$. Examples of such covariance kernels $\tilde k^i$ are given in Section \ref{applications}. For each of them, we will also provide the analytical expression of $k_0^i$.\\


Further, we also consider that the objective function can be rewritten as   
$$f(\xx) =f_0+ Z(\xx),$$
  where  $f_0$ is the constant mean of  $f(\xx)$ and $Z(\xx)$ is defined as (\ref{gp1}), where the covariance kernel $k(\xx,\txx)$ is given by (\ref{anovakernel}).
The choice of $k(\xx, \txx)$  is relevant here to propose a decomposition. Indeed, this definition implies the following properties:
\begin{enumerate}
\item If we set $Z_0 \sim \N{0}{\sigma^2}$, and, for all $u\in S$, $Z_u(\xx_u)\sim \GP(0,\sigma^2\prod_{i \in u} k_0^i(x_i, \tilde{x}_i))$ are independent processes, then, if

\begin{equation}\label{Zdecomposition}
Z(\xx) = Z_0 + \sum_{u \in S} Z_u(\xx_u), 
\end{equation}

we have that $Z(x)\sim \textrm{GP}(0,k(\xx,\txx))$. Indeed, as done in~\cite{durrande}, $k(\xx,\txx)$ can be decomposed as it follows,

\[
k(\xx,\txx)=\sigma^2+\sigma^2 \sum_{u\in S}\prod_{i\in u}k_0^i(x_i,\tilde x_i).
\]

\item Let us consider two sets  $u, v \in S$ such that $u \neq v$, and  $Z_u(\xx_u)\sim \GP(0,\sigma^2\prod_{i \in u} k_0^i(x_i, \tilde{x}_i))$,  $Z_v(\xx_v)\sim \GP(0,\sigma^2\prod_{i \in v} k_0^i(x_i, \tilde{x}_i))$. We have the following equalities  almost surely:
\begin{equation}\label{decZ}
\int Z_u(\xx_u)  \left( \prod_{i =1}^p p_{X_i}(x_i)  \right) \, d\xx  = 0, 
\end{equation}
and
\begin{equation}\label{orthoZ}
\int Z_u(\xx_u) Z_v(\xx_v)   \left( \prod_{i =1}^p p_{X_i}(x_i)  \right)  \, d\xx  = 0.
\end{equation}

Using Proposition \ref{prop1}, we know that the linear transformation $\int Z_u(\xx_u)  \left( \prod_{i =1}^p p_{X_i}(x_i)  \right) \, d\xx$ is Gaussian, and

\[
\int Z_u(\xx_u)  \left( \prod_{i =1}^p p_{X_i}(x_i)  \right) \, d\xx \sim \N{0}{\sigma^2 \int \prod_{i\in u}k_0^i(x_i,\tilde x_i) \prod_{i \in u} p_{X_i}(x_i) d \xx}.
\]

By replacing $k_0^i$ by its expression (\ref{k0i}), we deduce (\ref{decZ}). Now, for $u\neq v$, and $i\in u\setminus v$, we have,

\[
\int Z_u(\xx_u) Z_v(\xx_v)   \left( \prod_{i =1}^p p_{X_i}(x_i)  \right)  \, d\xx=\int Z_v(\xx_v)\left(  \int Z_u(\xx_u)   p_{X_i}(x_i)  dx_i \right) p_{\xx_{i^c}}    \, d\xx_{i^c},
\]

where $\xx_{i^c}$ is the complementary set of $x_i$, i.e. $\xx_{i^c}=(x_j)_{j\neq i}$. Again with Proposition \ref{prop1}, we conclude that (\ref{orthoZ}) is satisfied.
\end{enumerate}

\paragraph{Discussion about the choice of the covariance kernel when the input parameters are independent.} The  kernel  given in (\ref{anovakernel}) provides a relevant prior for  the sensitivity analysis when $P_\XX = \otimes_{i=1}^p P_{X_i}$. In this case, the Sobol index \cite{sobol} of $Z(\xx)$ --- modeling our prior knowledge about $f(\xx)$ ---  is given by:
\begin{displaymath}
S_u = \frac{V[Z_u(\XX_u)]}{V[Z(\XX)]},\quad \forall~u \in S,
\end{displaymath} 
where $Z_u$ checks Equalities (\ref{decZ}) and (\ref{orthoZ}).

By setting 
 $\sigma_u^2 = \sigma^2 \prod_{i \in u} \left( 1 +  k_0^i(x_i, x_i)  \right)$ and $\sigma_v^2 = \sigma^2 \prod_{i \in v} \left( 1 +  k_0^i(x_i, x_i)  \right)$ where $u \neq v \in S$,  we have:
\begin{displaymath}
\EZ \left[ V( Z_u( \XX_u ) \right] = \int \sigma_u^2 \prod_{i \in u} p_{x_i}(d x_i)
\end{displaymath}
and
\begin{displaymath}
\covZ \left(  V[ Z_u( \XX_u ) ] ,V[ Z_v( \XX_v ) ] \right) = 0.
\end{displaymath}

Therefore, we notice the sensitivity of $\XX_u$ in the model is monitored by $\sigma_u^2$, always strictly positive. This means that, through the decomposition (\ref{Zdecomposition}), we consider \emph{a priori} that every group of input variables is contributive in the model. We also consider \emph{a priori} that the sensitivity indices  are uncorrelated.


\section{ANOVA  decomposition of conditional Gaussian processes }\label{anovadec}

We propose in this section a representation of $f(\xx)$ as  a sum of increasing dimension Gaussian processes. Our main contribution is to consider the complete predictive distribution $[f(\xx) | \Zn =\y]$, given  by   (\ref{predictivedistribution}), and not only the predictive mean $\mu(\xx)$.
This allows us for quantifying the uncertainty due to the  meta-modeling  on the sensitivity indices estimation.

\subsection{ANOVA decomposition of the predictive mean}\label{ANOVAmean}
 This paragraph is dedicated to the functional decomposition of the predictive mean $\mu(\xx)$. Although it has been already developed and studied in~\cite{durrande}, we remind it here for the good understanding of the extension proposed in Paragraph \ref{ANOVAGP}.\\

 Remind that we consider the prior knowledge $Z(\xx) \sim \mathrm{GP}(0, k(\xx, \txx))$, where $k(\xx, \txx)$  is defined by Equations  (\ref{anovakernel})-(\ref{k0i}).
Following Paragraph \ref{GPRequations}, we consider the predictive distribution 
\begin{displaymath}
[f(\xx) | \Zn =\y] \sim \mathrm{GP}(\mu(\xx),s^2(\xx, \txx)), 
\end{displaymath}
where, from the definition of $k(\xx, \txx)$, we can decompose $\mu(\xx)$ as follows,
\begin{equation}\label{muxdecomposition}
\mu(\xx)= \mu_0 + \sum_{u \in S} \mu_u(\xx_u) , 
\end{equation} 
 with $\mu_0 = f_0 + \T \un \Kn^{-1}(\y -f_0\un)$  and
\begin{equation}\label{muxu}
\mu_u(\xx_u) = \prod_{i \in u}  
 \T \kk_{0,n}^i (\xx)  \Kn^{-1}(\y -f_0\un), \quad \forall~u\in S.
\end{equation}

 $\un$  the $n$-vector of $1$,  $\kn(\xx) =  [k (\xx, \xx^j)]_{j=1,\dots,n} $,  $\kk_{0,n}^i (x_i) = [k_{0}^i (x_i, x_i^j)]_{j=1,\dots,n}$ (see Equation (\ref{k0i})) and $\Kn =  [k (\xx^j, \xx^l)]_{j,l=1,\dots,n} $. 


Thanks to the property of the kernel $k_{0 }^i$, we can deduce that,
for all $u,v \in S$ with $u\neq v$:
\begin{displaymath}
\int \mu_u(\xx_u)    \left(\prod_{i=1}^p p_{X_i}(x_i)\right) d\xx  = 0,
\end{displaymath}
and
\begin{displaymath}
\int \mu_u(\xx_u)  \mu_v(\xx_v)  \left(\prod_{i=1}^p p_{X_i}(x_i)\right) d\xx  = 0.
\end{displaymath}

From this decomposition, \cite{durrande} deduce an analytical sensitivity measure, that we call $S_u^{\textrm D}$ here, to quantify the contribution of a given $\XX_u$ in the model. It is defined as,

\begin{equation}\label{indicedurrande}
S_u^{\textrm{D}}=\frac{V[\mu_u(\XX_u)]}{V[\mu(\XX)]},
\end{equation}

with  components $(\mu_u)_{u\in S}$ having the same properties as the ones of the Hoeffding expansion~\cite{hoeffding}. Thus, we can analyse the effect of each group of variables on the global variability, when the initial model $f$ is substituted to the predictive mean $\mu$.\\

However, when we only consider the predictive mean, we neglect an important part of information contained in the posterior variance. In addition, if the uncertainty of $\mu(\xx)$ is important, this means that the surrogate model does not ajust properly the objective function, leading to a wrong sensitivity analysis. \\
In the following part, we take into account the uncertainty of the meta-modeling by defining a functional ANOVA decomposition of a predictive distribution. Inspired by the work of~\cite{durrande}, we extend their work to a more general decomposition, and we also extend the sensitivity indices defined by (\ref{indicedurrande}) to the definition of sensitivity indices when input variables can be non independent.

\subsection{ANOVA decomposition of the conditional Gaussian processes}\label{ANOVAGP}

We saw in  the previous paragraph that the considered covariance kernel $k(\xx,\txx)$ (\ref{anovakernel}) leads to an ANOVA decomposition (\ref{muxdecomposition}) of $\mu(\xx)$  which is suitable to perform sensitivity analysis. However, as emphasized by \cite{oakley}, performing a sensitivity analysis based on a Gaussian process regression using only the predictive mean can be inappropriate. Indeed, in the framework of computer experiments, few observations $\y+f_0$ of $f(\xx)$ are available and thus the uncertainty on the meta-model $\mu(\xx)$ can be non-negligible.
For this reason, it is worth taking  into account the uncertainty on the meta-model and   to consider the complete predictive distribution. Nevertheless, it is non trivial to find the analogous of the ANOVA decomposition of the predictive mean to the predictive distribution. The proposition below allows for handling this issue~\cite{chiles1999geostatistics}.

\begin{prop}\label{anovaGP}
Let consider the random process $ f^n(\xx) $  defined as
\begin{equation}\label{Zetastar1}
f^n(\xx) = \mu(\xx) - {}^t \kn(\xx)\Kn^{-1}  \Zn +  Z(\xx),
\end{equation}
where $Z(\xx)\sim \text{GP}(0,k(\xx,\txx))$ with $k(\xx,\txx)=\sigma^2 \prod_{i =1}^p \left( 1 + k_0^i(x_i,\tilde{x_i}) \right)$, $\mu(\xx)$ is the predictive mean defined by (\ref{muxdecomposition}),  $\Zn = {}^t(Z(\xx^1), \dots, Z(\xx^n))$ is the Gaussian random vector corresponding to the value of $Z(x)$ at points in the experimental design set $ \D = \{\xx^1,\dots, \xx^n\}$,    $\kn(\xx) =  [k (\xx, \xx^j)]_{j=1,\dots,n} $ and $\Kn =  [k (\xx^i, \xx^j)]_{i,j=1,\dots,n} $.
Then, we have
\begin{equation}\label{memedistri}
f^n(\xx)  \sim [f(\xx) | \Zn =\y].
\end{equation}
\end{prop}

To get the proof of Proposition \ref{anovaGP}, the reader could refer to~\cite{chiles1999geostatistics}. 
Here, this result is of great interest since it allows for defining a Gaussian process $ f^n(\xx) $ distributed with respect to the predictive distribution $ [f(\xx) | \Zn =\y]$. Our goal is now to find a decomposition for $ f^n(\xx) $ which is suitable for performing sensitivity analysis.  
Following the  decompositions of $\mu(\xx)$ given in (\ref{muxdecomposition})-(\ref{muxu}) and of $Z(\xx)$ given in  (\ref{Zdecomposition}), the decomposition of $f^n$ is given in Proposition \ref{progpdecompo}.

\begin{prop}\label{progpdecompo}
Let $f^n(\xx)$ be the random process defined by (\ref{Zetastar1}) of Proposition \ref{anovaGP}. Then,

\begin{equation}\label{Zetastar}
 f^n(\xx)    =   f^n_0 + \sum_{u \in S}   f^n_u(\xx_u), 
\end{equation}
with 
\begin{equation}\label{decompZ}
\left\{\begin{array}{l}
 f^n_0  = \mu_0  -  \T \un \Kn^{-1} \Zn + Z_0,  \\
\\
 f^n_u(\xx_u) = \mu_u(\xx_u)  - \prod_{i \in u}  
 \T \kk_{0,n}^i (\xx)  \Kn^{-1} \Zn +  Z_u(\xx_u), \quad \forall u \in S.
\end{array}
\right.
\end{equation}
Furthermore, since  $k_{0 }^i$ defined in Equation (\ref{k0i})  is centered with respect to the marginal density $p_{X_i}$, the following properties  holds almost surely for all $u,v \in S$, $u \neq v$:
\begin{equation}\label{etanE}
\int  f^n_u(\xx_u)    \left(\prod_{i=1}^p p_{X_i}(x_i)\right) d\xx  = 0, 
\end{equation}
and
\begin{equation}\label{etanCOV}
\int   f^n_u(\xx_u)  f^n_v(\xx_v)   \left(\prod_{i=1}^p p_{X_i}(x_i)\right) d\xx  = 0.
\end{equation}

\end{prop}

The proof of proposition \ref{progpdecompo} is straightforward, by the decomposition of $\mu(\xx)$ given by (\ref{muxdecomposition})-(\ref{muxu}),  by the decomposition of $Z(\xx)$ given by (\ref{Zdecomposition}), and by the definition  of $k(\xx, \tilde \xx)$ given by (\ref{anovakernel}). The properties of the summands $(f_u^n)_u$ are also immediate. \\
%



\paragraph{Remark:} The suggested decomposition (\ref{Zetastar}) allows  for taking into account the meta-model uncertainty in the sensitivity index estimates. Furthermore, we highlight that it is easy to sample  with respect to the   distribution of $f^n (x)$. Indeed, from Proposition \ref{progpdecompo} we deduce that we can  perform it by sampling with respect to the distribution of   $Z(x) \sim \mathrm{GP}(0,k(\xx,\txx))$ and applying the linear transformation presented in (\ref{Zetastar1}). Moreover, to obtain a sample of $f^n_u(\xx_u) $ we just have to  sample   $Z_u(\xx_u) \sim \mathrm{GP}\left(0,  \prod_{i \in u} k_0^i(x_i, \tilde{x}_i)\right)$.

\section{Sensitivity measure definition for dependent   input variables}\label{si}

Now, we   adapt the methodology of \cite{li} to construct sensitivity indices for models with dependent inputs. First, we define a sensitivity measure based on the result of Proposition \ref{progpdecompo}. Then, we present an estimation procedure of the sensitivity measure in Paragraph \ref{indexestimation}. Finally, we present in Paragraph \ref{asymptoticnormality} how to take into account the uncertainty of the index estimation.

\subsection{Sensitivity measure definition}\label{indexdefinition}

As presented  in Section \ref{anovadec}, the  model function $f(\xx)$ is substituted to  $
  f^n(\xx)    =   f^n_0 + \sum_{u \in S}   f^n_u(\xx_u)$). Therefore, the global variance can now be decomposed as

\begin{equation}
V( f^n(\XX))=\sum_{u \in S}\bigg [V(  f^n_u(\XX_u))+\cov(  f^n_u(\XX_u), f^n_{u^c}(\XX_{u^c}))\bigg]
\end{equation}
where $   f^n_{u^c}(\XX)  =  f^n(\XX) - f^n_u(\XX_u)$.
 In this way, the sensitivity index associated to the group of variables $\XX_u$ is given by

\begin{equation}\label{indexmetamodel}
S_u^f=\frac{V[ f^n_u(\XX_u)]+\cov[  f^n_u(\XX_u),f^n_{u^c}(\XX )]}{V[ f^n(\XX)]}  
\end{equation}

We note that $S_u^f$ is defined on the probability space $(\Omega_Z, \mathcal{A}_Z,\mathbb P_Z)$ as $V$ and $\cov$ are the variance and covariance with respect to $P_\XX$. Therefore, $S_u^f$ integrates the uncertainty related to the meta-model approximation. In practice, the mode or the mean of the distribution $S_u^f$ may be used to get a scalar measure of sensitivity. \\

It also should be noted that this index is analogous with the Sobol one in the independent case when $f$ is replaced by $f^n$. Indeed, by  (\ref{etanE}) and (\ref{etanCOV}), we have

\begin{enumerate}
\item For $u\neq v\in S$, 
\begin{displaymath}
\cov \left( f^n_u(\XX_u), f^n_v(\XX_v) \right) = 0.
\end{displaymath}

\item For a given $u\in S$, by integrating $f^n$ with respect to the distribution of the all inputs except the ones indexed by $u$, we have that
\begin{displaymath}
f^n_u(\XX_u) = \E \left[f^n(\XX) | \XX_u \right] +\sum_{\substack{v \subset  u \\ v \neq u}} (-1)^{|u| - |v|}\E \left[f^n(\XX) | \XX_v \right].
\end{displaymath}
\end{enumerate}
 Hence, when $P_\XX=P_{X_1}\otimes \cdots \otimes P_{X_p}$, 
 \[
 S_u^f =\frac{V( \E \left[f^n(\XX) | \XX_u \right])+\sum_{\substack{v \subset  u \\ v \neq u}} (-1)^{|u| - |v|} V(\E \left[f^n(\XX) | \XX_v \right])}{V[ f^n(\XX)]}.
 \]

For models with dependent inputs, the parametric  functional ANOVA decomposition mentioned in Section \ref{introduction}~\cite{li,caniouth,chastaing2}  requires to estimate $2^p$   components, which is hardly achievable in practice when $p$ gets large. Thus, their decomposition must be truncated but in this case we loose a part of model information.

Here, we have accessed to the approximation $f^n(\xx)$  of $f(\xx)$ without processing all the terms $f_u^n(\xx_u)$, $u \in S$,  thanks to the equality $ f^n(\xx) = \mu(\xx) - {}^t \kn(\xx)\Kn^{-1}  \Zn +  Z(\xx)$
where $Z(\xx)\sim \text{GP}(0,k(\xx,\txx))$, and $ \mu(\xx)  =  f_0 + {}^t \kn(\xx)\Kn^{-1}(\y -f_0\un)$. Furthermore, for any  $u \in S$, we have an explicit expression of $f_u^n(\xx_u)$ given by (\ref{decompZ}). Thus, as $f^n=f_u^n+f_{u^c}^n$, the complementary summand $f_{u^c}^n$ of $f_{u}^n$ can be easily deduced, avoiding the truncation error in the estimation.




\subsection{Estimation procedure}\label{indexestimation}

The aim of this section is to provide an efficient numerical estimation of the sensitivity measure $S_u^f$, for a given $u\in S$. As already mentioned,  $S_u^f$ lies in $(\Omega_Z,\mathcal A_Z,\mathbb P_Z)$. Further, we describe a numerical procedure for only one realization $S_u^f$. In practice, this procedure is repeated $N_s$ times to take into account the uncertainty of the meta-model.\\

We empirically estimate  the variance and the covariance  involved in (\ref{indexmetamodel}) with a   Monte-Carlo integration.  Therefore, we consider the following estimator from $m$ realizations $\TT = (\tx_1^j,\cdots,\tx_p^j)_{j=1,\dots,m}$ of the  random variable $\XX$ defined on the probability space $(\Omega_\XX, \mathcal{A}_\XX, \mathbb{P}_\XX)$:
\begin{equation}\label{estimatorMC}
S_{u,m}^f= \frac{
\frac{1}{m}\sum_{j=1}^m  f^n_u(\tx^j_u) ^2 - \left( \bar{f}^n_u \right)^2 + \frac{1}{m}   \sum_{j=1}^m   f^n_u(\tx^j_u)  f^n_{u^c}(\tx^j)  - \bar{f}^n_u\bar{f}^n_{u^c}
}{
\frac{1}{m}\sum_{j=1}^m  f^n(\tx^j ))^2 - \left( \frac{1}{m}\sum_{j=1}^m  f^n(\tx^j ) \right)^2
}
\end{equation}
where $\bar{f}^n_u = \frac{1}{m}\sum_{j=1}^m  f^n_u(\tx^j_u) $,  $\bar{f}^n_{u^c} = \frac{1}{m}\sum_{j=1}^m  f^n_{u^c}(\tx^j) $ and $\tx^j_u = (\tx^j_i)_{i \in u}$, $j = 1,\dots,m$. We point out  that  $S_{u,m}^f$  lies in the product probability space $(\Omega_Z, \mathcal{A}_Z, \mathbb P_Z)$.
Therefore, we generate several realizations of $S_{u,m}^f$ to get an estimate of it.
This procedure is described further below.\\

First, let us denote  $\MM = \begin{pmatrix} \TT \\ \D  \end{pmatrix}$, where $\D = (\xx ^j)_{j=1,\dots,n}$ is the experimental design set. For $u \in S$, we denote $\MM_u = \begin{pmatrix} \TT_u \\ \D_u  \end{pmatrix}$, where $\TT_u =  (\tx_u ^j)_{j =1,\dots,m}$ and $\D_u =  (\xx ^j_u)_{j =1,\dots,n}$.  

\begin{enumerate}
\item \label{step1} To get a realization of $f_u^n(\xx_u)$, we generate a sample from the distribution of 
\begin{displaymath} 
Z_u(\xx_u) \sim \mathrm{GP}\left(0,\prod_{i \in u}k_0^i(x_i,\tilde{x}_i) \right), 
\end{displaymath}  
on   $\MM$ with the following procedure:
\begin{enumerate}
\item Compute   $\K_{0,m}^u=\bigodot_{i\in u}k_0^i\left( \MM_i, \MM_i \right)$ and the Cholesky decomposition $\LL_{0,m}^u$ of $\K_{0,m}^u$ where    $\bigodot$ stands for the term-wise matrix product.    $\K_{0,m}^u$  is the covariance matrix  of $Z_u(\xx_u)$ at points in $\MM_u$.
\item Generate one realization $z_u(\MM_u)$ of $ Z_u(\xx_u)$ on $\MM_u$ from the Cholesky decomposition of $\K_{0,m}^u$ (see \cite{rasmussen} Appendix A.2) with,

\begin{equation}\label{samplezu}
z_u(\MM_u) = \LL_{0,m}^u \eepsilon_u,
\end{equation}
where $\eepsilon_u $ is a sample generated from the distribution $   \N{0}{I_{n+m}}$ where $I$ is the identity matrix of size $(n+m) \times (n+m)$. 
\end{enumerate}
\item \label{step2} To generate a realization of $f_{u^c}(\xx)$ on   $\MM$, we generate a sample from the distribution of 
\begin{displaymath} 
Z_{u^c}(\xx) \sim \mathrm{GP}\left(0,k(\xx, \txx) - \prod_{i \in u}k_0^i(x_i,\tilde{x}_i) \right), 
\end{displaymath}
on   $\MM$, with the following procedure:
\begin{enumerate}
\item Compute  $\K_{m}= k \left( \MM, \MM \right)$ and the Cholesky decomposition  $\LL_{0,m}^{u^c}$ of $\K_{m} - \K_{0,m}^u$. $k \left( \MM, \MM \right)$ is the covariance matrix of $Z(\xx)$  at points in $\MM$. 
\item Generate one realization $z_{u^c}(\MM)$ of $ Z_{u^c}(\xx)$ on $\MM$ with
\begin{equation}\label{samplezuc}
z_{u^c}(\MM) = \LL_{0,m}^{u^c}  \eepsilon_{u^c},
\end{equation}
where $\eepsilon_{u^c}$ is sampled from the distribution $\N{0}{I_{n+m}}$.
\end{enumerate}
\item We  deduce a  sample $z(\MM)$  of $ Z(\xx)$ on $\MM$ with:
\begin{equation}
z(\MM) = z_u(\MM_u) + z_{u^c}(\MM) .
\end{equation}
Moreover,  as $M=\begin{pmatrix} \TT \\ \D \end{pmatrix}$ and $M_u=\begin{pmatrix} \TT_u \\ \D_u \end{pmatrix}$, $z(\MM) $, $z_u(\MM_u)$ and $z_{u^c}(\MM) $ can be rewritten in the following forms:
\begin{displaymath}
z(\MM) = \begin{pmatrix}z( \TT) \\ z(\D)  \end{pmatrix},  \, \, 
z_u(\MM_u) = \begin{pmatrix}z_u( \TT_u) \\ z_u(\D_u)  \end{pmatrix}  \,\,  \mathrm{and} \, \,
z_{u^c}(\MM) = \begin{pmatrix}z_{u^c}( \TT) \\ z_{u^c}(\D)  \end{pmatrix}.
\end{displaymath}
\item We can deduce the samples  $\tilde f^n(\TT)$,   $\tilde f^n_u(\TT_u)$ and $\tilde f^n_{u^c} (\TT)$   of  respectively  $ f^n(\xx)$,  $ f^n_u(\xx_u)$ and $ f^n_{u^c} (\xx)$   on $\TT$ with the following formulas:
\begin{equation}\label{decompZ2}
\left\{\begin{array}{l}
\tilde f^n(\TT)  = \mu(\TT) -   k(\TT, \D)\Kn^{-1}  z(\D) +  z(\TT),  \\
\tilde  f^n_u(\TT_u) = \mu_u(\TT_u)  - 
   \left( \bigodot_{i\in u}k_0^i\left( \TT_i, \D_i \right) \right)  \Kn^{-1} z(\D)  + z_u(\TT_u) ,  \\
\tilde f^n_{u^c} (\TT) = \tilde f^n(\TT)  -  \tilde f^n_u(\TT) 
\end{array}
\right.
\end{equation}
where $\Kn = k(\D,\D)$, and

\[
\left\{
\begin{array}{l}
\mu(\TT)={}^t k(\TT,\D) \Kn^{-1} (\y-f_0\un)+f_0\\
\mu_u(\TT_u)= \bigodot_{i\in u}{}^t k_0^i\left( \TT_i, \D_i \right) \Kn^{-1}(\y-f_0\un)
\end{array}
\right.
\] 
We note that $k(\TT, \D)$ and $ \left( \bigodot_{i\in u}k_0^i\left( \TT_i, \D_i \right) \right) $ have   been already computed with $\K_{m}$ and $\K_{0,m}^u$ (Step \ref{step1} and \ref{step2}) as,

$$
\K_m = k(\MM, \MM) = \begin{pmatrix} k(\TT, \TT)  &  k(\TT, \D) \\ k(\D, \TT)  & k(\D, \D) \end{pmatrix},
$$  
and  
$$\K_{0,m}^u=\bigodot_{i\in u}k_0^i\left( \MM_i, \MM_i \right) = \begin{pmatrix} \bigodot_{i\in u}k_0^i\left( \TT_i, \TT_i \right) & \bigodot_{i\in u}k_0^i\left( \TT_i, \D_i \right) \\ \bigodot_{i\in u}k_0^i\left( \D_i, \TT_i \right) & \bigodot_{i\in u}k_0^i\left( \D_i, \D_i \right)\end{pmatrix}.
$$

\item  We deduce that a sample $s_{u,m}^f$  of $S_{u,m}^f$ is given by

\begin{equation}\label{Ssample}
s_{u,m}^f= \frac{
\frac{1}{m}\sum_{j=1}^m  \tilde f^n_u(\tx^j_u) ^2 - \left( \bar{\tilde{f}}^n_u \right)^2 + \frac{1}{m}   \sum_{j=1}^m \tilde  f^n_u(\tx^j_u) \tilde f^n_{u^c}(\tx^j)  - \bar{\tilde{f}}^n_u \bar{\tilde{f}}^n_{u^c}
}{
\frac{1}{m}\sum_{j=1}^m  \tilde f^n(\tx^j ))^2 - \left( \frac{1}{m}\sum_{j=1}^m  \tilde f^n(\tx^j ) \right)^2
}, 
\end{equation}
where $\bar{\tilde{f}}^n_u = \frac{1}{m}\sum_{j=1}^m \tilde f^n_u(\tx^j_u) $ and  $\bar{\tilde{f}}^n_{u^c} = \frac{1}{m}\sum_{j=1}^m \tilde  f^n_{u^c}(\tx^j) $.
\end{enumerate}

\subsection{Asymptotic normality of the sensitivity index estimator}\label{asymptoticnormality}

We have presented  in the previous paragraph  a procedure to sample $S_{u,m}^f$ defined in Equation (\ref{estimatorMC}). However, for given realizations of random processes $  f^n_u(\XX_u)$ and $ f^n_{u^c}(\XX )$, the estimated sensitivity measure comes from a Monte-Carlo integration and thus integrates a Monte-Carlo error. The purpose of this paragraph is to quantify it. A natural approach is to use an asymptotic normality result as stated in  Proposition \ref{normalityasymptotic}.

\begin{prop}\label{normalityasymptotic}
For $u\in S$, let us consider  the respective  realizations of  $f^n_u(\XX_u)$, $f^n_{u^c}(\XX )$ and $f^n(\XX)$ denoted by $\tilde  f^n_u(\XX_u)$, $\tilde f^n_{u^c}(\XX )$ and $\tilde f^n(\XX)$ respectively. We denote the theoretical   sensitivity measure for $\XX_u$ associated to $\tilde  f^n$ by
\begin{displaymath} 
 s_u^f=\frac{V(  \tilde  f^n_u(\XX_u))+\cov(  \tilde  f^n_u(\XX_u), \tilde f^n_{u^c}(\XX ))}{V( \tilde f^n(\XX))} , 
\end{displaymath}
 where we assume that $V(\tilde f^n(\XX)) \neq 0$. Suppose also  that $\E \left[ \tilde f_u^n(\XX)^4\right] < \infty$ for all $u \in S$. Then, for any $u \in S$, we have, when $m \rightarrow \infty$:
\begin{equation}\label{tcl}
\sqrt{m} \left(  s_{u,m}^f - s_u^f \right) \stackrel{\mathcal{L}}{\longrightarrow}  \mathcal{N} \left( 0 , \T(\bigtriangledown \pphi(\mmu) )\GGamma \bigtriangledown \pphi(\mmu)  \right)
\end{equation}
where $\mmu = \E(\U)$, $\GGamma = V(\U)$, 
\begin{displaymath}
\U = \begin{pmatrix}
\tilde f^n_u(\XX_u)  &
\tilde f^n_{u^c}(\XX)  &
\tilde f^n(\XX)  &
\tilde f^n_u(\XX_u) ^2 &
\tilde f^n(\XX) ^2  &
\tilde f^n_u(\XX_u)\tilde f^n_{u^c}(\XX)
\end{pmatrix},
\end{displaymath}
and
\begin{displaymath}
 \pphi(u_1, u_2, u_3, u_4, u_5, u_6) = \frac{u_4 - u_1^2 + u_6 - u_1 u_2}{u_5 - u_3^2}.
\end{displaymath}
\end{prop}
The proof of Proposition \ref{normalityasymptotic} is straightforward by using the Delta method  (see \cite{vandervaart}). We highlight that the terms $\mmu$ and $\GGamma $ in Proposition \ref{normalityasymptotic} can be estimated from the sample $\TT = (\tx^j)_{j=1,\dots,m}$ used in the Monte-Carlo integration (\ref{estimatorMC}). \\



In practice, we use the asymptotic result given in Proposition \ref{normalityasymptotic} to estimate the Monte Carlo error. Thus, to take into account both the uncertainty of the surrogate model and the one of the Monte Carlo integration, we proceed as follows,

\begin{enumerate}
\item \label {item1} Generate $\tilde f^n(\TT)$,   $\tilde f^n_u(\TT_u)$ and $\tilde f^n_{u^c} (\TT)$  from the sample $\TT$ with the estimation procedure given in Paragraph \ref{indexestimation}.
\item \label{item2} Generate a sample of size $K$ from the limit distribution $ \mathcal{N} \left( 0 , \T(\bigtriangledown \pphi(\hat \mmu) ) \hat \GGamma \bigtriangledown \pphi(\hat \mmu)  \right)
$ where $\hat \mmu = \frac 1 m \sum_{j=1}^m \U(\tx_u^j)$, $\hat \GGamma = \frac 1 m \sum_{j=1}^m [\U(\tx_u^j)-\hat \mmu]^2$, and 
\begin{displaymath}
\U = \begin{pmatrix}
\tilde f^n_u  &
\tilde f^n_{u^c}  &
\tilde f^n  &
(\tilde f^n_u)^2 &
(\tilde f^n) ^2  &
\tilde f^n_u \tilde f^n_{u^c}
\end{pmatrix}.
\end{displaymath}
Thus, a Monte Carlo error is obtained for one realization $ s_{u,m}^f$.
\item Repeat Steps \ref{item1}-\ref{item2} $N_s$ times to take into account the uncertainty of the surrogate model.
\end{enumerate}



\section{Applications}\label{applications}

We illustrate in this section our sensitivity measure on academic and industrial examples. First, we present explicit examples of covariance kernels $k_0^i$, and a procedure to estimate their parameters. Further, we illustrate the estimation procedure of Paragraph \ref{indexestimation} through several numerical applications.

\subsection{Example of covariance kernels}\label{Examplek0}

Here, we analytically compute  zero mean kernels for two usual kernels associated with uniform distributions. First, let us consider that (see Equation (\ref{k0i})):
\begin{displaymath} 
k_0^i(x, \tilde x) = \tilde{k}^i(x, \tilde x) - \frac{\int {\tilde{k}^i(x, u) p_{X_i}(u)du } \int{ \tilde{k}^i(v, \tilde x) p_{X_i}(v)dv} }{\int  \int \tilde{k}^i(u, v) p_{X_i}(u) p_{X_i}(v) du   dv  }, \quad x, \tilde x \in \R.
\end{displaymath}

\paragraph{Example 1:} We consider an exponential kernel  for $\tilde{k}^i(x, \tilde x)$ with an uniform marginal $p_{X_i}$, namely,
\begin{displaymath}
k^i(x, \tilde x) =  \exp \left(-\frac{1}{2} \frac{ |x-\tilde x|}{\theta_i} \right), \quad \theta_i > 0, 
\end{displaymath}
and 
\begin{displaymath}
p_{X_i} \sim \mathcal{U}\left(a_{i}, b_{i}\right).
\end{displaymath}
Then, the corresponding covariance kernel $k_0^i(x, \tilde x) $ is given by:
\begin{displaymath}
\begin{array}{ll}
k^i_0(x, \tilde x)=&  \exp \left(-\frac{1}{2} \frac{ |x-\tilde x|}{\theta_i} \right)- \frac{\theta_i}{b_i-a_i+2\theta_i \left( \exp(-\frac{1}{2}\frac {b_i-a_i} {\theta_i})-1 \right)}\times\\
& \left(2-\exp(-\frac{1}{2}\frac {x-a_i} {\theta_i})-\exp(-\frac{1}{2} \frac {b_i-x} {\theta_i})\right)\cdot \left(2-\exp(-\frac{1}{2}\frac {\tilde x-a_i} {\theta_i})-\exp(-\frac{1}{2} \frac {b_i-\tilde x} {\theta_i})\right)
\end{array}
\end{displaymath}
We note that the exponential kernel is stationary --- i.e. it is invariant under translations in the input parameter space --- and corresponds to the covariance of an Ornstein-Uhlenbeck process. Furthermore, the corresponding process is continuous in mean square sense and nowhere differentiable. Therefore this kernel is appropriate for rough function $f(\xx)$.
\paragraph{Example 2:} We consider a Gaussian  kernel  for $\tilde{k}^i(x, \tilde x)$ with an uniform marginal $p_{X_i}$, namely,
\begin{displaymath}
k^i(x, \tilde x) = \exp \left(-\frac{1}{2} \frac{ (x-\tilde x)^2}{\theta_i^2} \right), \quad \theta_i > 0, 
\end{displaymath}
and 
\begin{displaymath}
p_{X_i} \sim \mathcal{U}\left(a_{i}, b_{i}\right).
\end{displaymath}
Then, the corresponding covariance kernel $k_0^i(x, \tilde x) $ is given by:
\begin{displaymath}
\begin{array}{ll}
k^i_0(x, \tilde x)=&  \exp \left(-\frac{1}{2} \frac{ (x-\tilde x)^2}{\theta_i^2} \right) - A(x) A(\tilde x) / B, 
\end{array}
\end{displaymath}
where 
\begin{displaymath}
A(x) = -\frac{\sqrt{\pi}}{\sqrt{2}} \theta_i \mathrm{erf}\left(    \frac{a_i -x}{ \theta_i{\sqrt{2}}}\right)
      + \frac{\sqrt{\pi}}{\sqrt{2}}  \theta_i \mathrm{erf}\left(  \frac{b_i-x}{\theta_i \sqrt{2}}\right), 
\end{displaymath}
\begin{displaymath}
B =-2\theta_i^2+\theta_i \sqrt{2} \mathrm{erf}\left (\frac{a_i-b_i}{\theta_i \sqrt{2}} \right)\sqrt{\pi}(a_i-b_i)
      +2 \exp \left(-\frac{1}{2} \frac{ \left(a_i -b_i\right)^2}{\theta_i^2}\right)\theta_i^2,
\end{displaymath}
and the error function is given by
\begin{displaymath}
 \mathrm{erf}(x)  = \frac{2}{\sqrt{ \pi}} \int_0^x \exp(-t^2)dt.
\end{displaymath}

We note that the Gaussian kernel corresponds to   processes   infinitely continuously differentiable in  mean square sense. Therefore this kernel is appropriate for very smooth function $f(\xx)$.

Closed form expressions can also be derived  for $5/2$-Mat\'ern and $3/2$-Mat\'ern covariance kernels (see \cite{S99}) by straightforward calculations.  Due to their complex expression, they are not presented here. Though, we note that these kernels correspond respectively to once and twice continuously differentiable processes in mean square sense. Therefore, they could be a relevant compromise between the exponential and the Gaussian kernels.

\subsection{Covariance kernel  parameter estimation}\label{parameterestimation}

We deal in this section with the estimation of the model parameters using a maximum likelihood method. Let us consider the covariance kernel
\begin{displaymath} 
k(\xx, \txx) = \sigma^2 \prod_{i=1}^p \left( 1 +  k_0^i(x_i, \tilde{x}_i)  \right), 
\end{displaymath}

 where $k_0^i(x_i, \tilde{x}_i)$ is one  of the covariance kernels given in  Example 1 or Example 2 of Paragraph \ref{Examplek0}. Therefore, the parameters to be estimated are the variance parameter $\sigma^2$, the mean $f_0$ and the hyper-parameter $\boldsymbol{\theta} = ( \theta_i)_{i=1,\dots,p}$ of $( k_0^i(x_i, \tilde{x}_i))_{i=1,\dots,p}$.\\

First, let us consider the maximum likelihood estimate of $f_0$:
\begin{equation}
\hat f_0 = \left( {}^t \un \Kn^{-1} \un \right)^{-1}{}^t \un \Kn^{-1} \y, 
\end{equation}
where $\y= (f(\xx^i))_{i=1,\dots,n} $ and $\Kn = [k(\xx^i, \xx^j)]_{i,j=1,\dots,n}$ is the covariance matrix of the observations at points $\D = (\xx^i)_{i=1,\dots,n}$, with $\xx^i \in \R$, for all $i=1,\dots,n$. Then, we substitute $\hat f_0$ in the likelihood and maximize it with respect to $\sigma^2$. We obtain the following maximum likelihood estimate of $\sigma^2$:
\begin{equation}
\hat \sigma ^2 = \frac{{}^t (\y -\hat f_0 \un )  \Kn^{-1}  (\y -\hat f_0 \un )   }{n}.
\end{equation}
Finally, we substitute $\sigma^2$ with $\hat \sigma ^2$ in the likelihood to obtain the   marginal likelihood:
\begin{equation}\label{marginallikelihood}
\mathcal{L}(\boldsymbol{\theta} ; \y) = n \log (\hat \sigma ^2 ) + \log(\det \Kn  ).
\end{equation}
The estimate $\hat{ \boldsymbol{\theta}}$ of $\boldsymbol{\theta}$ is obtained by minimizing (\ref{marginallikelihood}) with respect to $\boldsymbol{\theta}$.  In practice, we use an evolutionary algorithm coupled with a BFGS (Broyden-Fletcher-Goldfarb-Shanno) procedure (see \cite{avriel2003nonlinear}).

\subsection{Academic example: the Ishigami function}

 Let us consider the Ishigami function:
\begin{equation}
z(x_1, x_2, x_3) = \sin (x_1) + 7 \sin(x_2)^2 + 0.1 x_3^4 \sin(x_1)
\end{equation}
with $(x_1, x_2, x_3) \in [-\pi, \pi]^3$. This function is  a classical tabulated function for sensitivity analysis~\cite{saltelli}. 

\subsubsection{Gaussian process regression model building}

First of all,  let us present the meta-model building. The considered experimental design set is a Latin-Hypercube-Sampling (LHS)  \cite{St87} of $n=150$ points optimized with respect to the maximin criterion. This criterion maximizes the minimum distance between the points.
We consider the Gaussian covariance kernel presented as Example 2 of Paragraph  \ref{Examplek0}.
 The maximum likelihood estimates of the model parameters are given below (see Paragraph \ref{parameterestimation}):
\begin{displaymath}
\hat \ttheta  =\begin{pmatrix}1.98 &   1.44  &   1.63 \end{pmatrix}, ~~\hat{\sigma}^2 = 16.50,~~\hat{f}_0 = 3.40.
\end{displaymath}

The efficiency of the model is assessed on a test set $\XX_\mathrm{test}$ of size $n_t$ uniformly spread on $[-\pi, \pi]^3$ with the following  coefficient:
\begin{displaymath}
Q^2 =  1 - \frac{
\sum_{\xx \in \XX_\mathrm{test}} \left( \mu(\xx) - f(\xx) \right)^2
}{
\sum_{\xx \in \XX_\mathrm{test}} \left( \mu(\xx) - \bar{f} \right)^2
}, \quad \bar{f} = \frac{\sum_{\xx \in \XX_\mathrm{test}} f(\xx)}{n_t},
\end{displaymath}

where  $\mu(\xx)$ is the predictive mean given in (\ref{muxdecomposition}). The coefficient $Q^2 $ represent the part of the  model discrepancy explained by the meta-model. The closer to 1, the more efficient is  the meta-model.
 Here, the estimated efficiency is $Q^2 = 98.2 \%$. Then, we have an accurate  meta-model.

\subsubsection{Ishigami function with independent inputs}

We consider the product measure $P_\XX = P_{X_1} \otimes P_{X_2} \otimes P_{X_3}$ with $X_i$ uniformly  distributed on $[-\pi,\pi]$, i.e. $X_i \sim \mathcal{U}(-\pi, \pi)$, for $i=1,2,3$. In this case, the theoretical  sensitivity indices coincide  to the classical Sobol indices (see \cite{sobol}). Their values are indicated in Table \ref{tab_ishigami_independent}. The purpose of this   paragraph is to study the relevance of the suggested indices in the case of independent inputs. 

To perform the Monte-Carlo integration presented in Paragraph \ref{indexestimation}, we  generate a  sample $\TT = (\tx^j)_{j=1,\dots,m}$  of  $m = 10,000$ points with respect to the product measure $P_\XX$. Further, we generate $N_s = 200$ realizations  (see $s^f_{u,m}$ in  Equation  (\ref{Ssample})) of  the estimator $S^f_{u,m}$ of the sensitivity measure $S ^f_u$ with $u \in \{\{1\}, \{2\}, \{3\}, \{1,2\}, \{1,3\}, \{2,3\}\}$. \\

To get a scalar quantity estimate $\hat{S} ^f_u$  of  $S ^f_u$, we consider that $\hat{S} ^f_u$ is the mode of the probability density estimate  of  the $N_s = 200$ realizations of  $S^f_{u,m}$. We note that the estimate of the probability density  is based on a normal kernel function with a window parameter  
optimal for estimating a normal density  \cite{Bowman}.

The estimated indices are given in Table \ref{tab_ishigami_independent}. Furthermore, we provide the confidence intervals of each estimators using the procedure given in Paragraph \ref{asymptoticnormality} that allows to take into account the meta modeling and the Monte Carlo errors. To do that, we consider a sample of size $K=200$ for each realization.

\begin{table}[H]
\begin{center}
\begin{tabular}{|c|c|c|c|c|c|c|}
\hline
Index & $S_1$ & $ S_2$  & $ S_3$ & $ S_{12}$ & $ S_{13}$ & $ S_{23}$ \\
\hline
Analytical & 0.314  &0.442   &  0 & 0  & 0.244 & 0 \\
\hline\hline
Estimate &  0.310 & 0.447 &  0.000 & 0.001  & 0.238 & 0.000 \\
\hline
2.5\%-quantile  & 0.308    & 0.430 &    -0.001    & -0.001 & 0.231  & -0.001    \\
\hline
97.5\%-quantile & 0.320  & 0.452 &  0.001 & 0.002 &  0.253   & 0.001 \\
\hline
\end{tabular}
\end{center}
\caption{Sensitivity measure estimates for the Ishigami function with independent input parameters.}
\label{tab_ishigami_independent}
\end{table}
We see that the estimated sensitivity measures fit the theoretical ones. This emphasizes the efficiency of the suggested estimation procedure. 

To show the relevance of the estimated 95\%-confidence intervals, we reiterate the presented procedure with 500 Gaussian process regression models built from different maximin-LHS design sets. For each design set, the parameters $\ttheta$, $\sigma^2$ and $f_0$ are estimated with a maximum likelihood method. We thus have 500 estimated 95\% confidence intervals and we verify whether they include the true index or not. The ratio of intervals including the true indices, also called the coverage rate, has to be close to 95\%.
The results of this procedure is presented in Table \ref{tab_ishigami_coverage}. Moreover, 
these intervals are compared with the ones considering only the meta-modeling error (i.e. without using the procedure presented in Paragraph \ref{asymptoticnormality} to evaluate the Monte-Carlo integration error). \\
Furthermore, to show the issue involving the meta-modeling error, we compare the coverage rate of the empirical estimation $\hat S_u^{\textrm{D}}$ of the sensitivity index $S_u^{\textrm{D}}$ proposed by Durrande \etal(see  Paragraph \ref{ANOVAmean})  to the one of our sensitivity measure.  To evaluate the Monte-Carlo error of $\hat S_u^{\textrm{D}}$, we  apply the procedure presented in Paragraph \ref{asymptoticnormality}. The results are presented in Table \ref{tab_ishigami_coverage}. 

%

\begin{table}[!ht]
\begin{center}
\begin{tabular}{|c|c|c|c|c|c|c|c|}
\hline
Index & Error & $ {S}_1$ & $  S_2$  & $  S_3$ & $  S_{12}$ & $  S_{13}$ & $  S_{23}$ \\
\hline
\multirow{2}{1cm}{$\hat S_u^f$} &
MC + Meta-model & 0.95  &  0.97  &  0.91   &  0.76   &  0.81 &   0.74   \\
\cline{2-8}
 & Meta-model &   0.87 &   0.94  &  0.87  &  0.08  &  0.29   &  0.21 \\
\hline\hline
\multirow{1}{1cm}{$\hat S_u^D$} &
MC + predictive mean &    0.67 &   0.65   &  0.39  & 0.00   & 0.26   &  0.05  \\
\hline
\end{tabular}
\end{center}
\caption{Coverage rates of 500 empirical confidence intervals. The theoretical confidence interval  is 95\%. The coverage rates for the suggested confidence intervals taking into account the uncertainty of both the meta-model approximation and the Monte-Carlo integration are labeled ``MC+meta-model'' ; the ones taking into account only the meta-model error are labeled ``meta-model'' ; the ones taking into account only the Monte-Carlo error and using the predictive mean are labeled ``MC+predictive mean''.}
\label{tab_ishigami_coverage}
\end{table}

We see in Table \ref{tab_ishigami_coverage} that the   empirical confidence intervals obtained  with the suggested procedure are better than those which only take into account the meta-model   or the Monte-Carlo error. In particular, the confidence intervals found with the estimator $\hat S_u^\mathrm{D}$   and considering only the Monte-Carlo error are widely underestimated. However, we see that they are all  underestimated   for the second order indices. Especially for the   indices  corresponding to  the non-influent interactions, namely, $S_{12}$ and $  S_{23}$. The underestimation could be due to a poor learning of the interactions by the meta-model.

\subsubsection{Ishigami function with perfectly  correlated inputs}

We present here a sensitivity analysis   with $P_{X_i} \sim \mathcal{U}(-\pi, \pi)$, $i=1,2,3$,  and where we assume that $X_1 = X_2$ and $X_1$, $X_2$   independent of $X_3$.\\
Therefore,  we can either perform a sensitivity analysis considering only two independent variables (Case i) or perform a sensitivity analysis with three input variables where two  of them  have a perfect positive linear relationship (Case ii). We use the classical Sobol indices for Case i, and we perform our procedure for Case ii.

As the two sensitivity analyses formally correspond to the same underlying function, it should have a connection between them. The purpose of this paragraph is to  numerically observe it.
 Since $X_1 = X_2$, we can consider the following function:
\begin{displaymath}
z_{(i)}^{\mathrm{sob}}(X_1, X_3) = \sin (X_1) + 7 \sin(X_1)^2 + 0.1 X_3^4 \sin(X_1),
\end{displaymath}
with $X_1$ independent of $X_3$. Further, we denote  $\hat S_{u}^{\mathrm{sob}}$, for $u\in \{ \{1\}, \{3\}, \{13\} \}$ the estimators of the Sobol index. Our indices are denoted $\hat S_u^f$, for $u \in \{ \{1\}, \{2\},\{3\}, \{12\}, \{13\},\{23\}\}$, as we consider the mode of the distribution estimate. Results are given in Table \ref{tab_ishigami_dependent_X1X2}.

%
\begin{table}[!ht]
\begin{center}
\begin{tabular}{|c|c|c|c|c|c|c|}
\hline
Index & $\hat{S}_1$ & $\hat S_2$  & $\hat S_3$ & $\hat S_{12}$ & $\hat S_{13}$ & $\hat S_{23}$ \\
\hline
$\hat S_u^f $ &  0.308  & 0.439  & 0.001   & 0.001   &  0.238  &   0.012 \\
\hline
$\hat S_{u}^{\mathrm{sob}}$ & 0.751 & - & 0.001& - & 0.245 & -\\
\hline
\end{tabular}
\end{center}
\caption{Sensitivity measure estimates for the Ishigami function with $X_1 = X_2$, $X_1$, $X_2$   independent of $X_3$ and $P_{X_i} \sim \mathcal{U}(-\pi, \pi)$, $i=1,2,3$.}
\label{tab_ishigami_dependent_X1X2}
\end{table}
We see in  Table \ref{tab_ishigami_dependent_X1X2} that we empirically found that $\hat S_{1}^{\mathrm{sob}} \approx \hat{S}_1^f + \hat S_2^f + \hat S_{12}^f$, $\hat S_{3}^{\mathrm{sob}}  \approx \hat S_3^f$ and $\hat S_{13}^{\mathrm{sob}} \approx \hat S_{13}^f + \hat S_{23}^f$. Therefore, we numerically observe a direct correspondence between the classical sensitivity analysis for independent inputs and the suggested one for dependent inputs. 

This connection strengthen the relevance of the considered index since, in the independent case, the Sobol indices are  commonly accepted as a good  measure of sensitivity.
However, the connection is only established when the considered model  can  be reduced to an equivalent model which have  independent inputs. 
For general cases, the interpretation will be much more complex (see Paragraph \ref{industrial}).

\subsubsection{Modeling dependence with copulas}

To define the dependence among random variables, it is usual to use the copula functions~\cite{nelsen}. Indeed, a copula function aims to join the joint distribution of a set of variables to its marginal distributions. If the cumulative distribution function (cdf) of $\XX$ is denoted $F_\XX$, and $F_1,\cdots,F_p$ are the respective marginal cdf of $X_1,\cdots,X_p$, then there exists a copula $p$-dimensional $C$ such that, for all $\xx\in \R^p$,

\[
F_\XX(\xx)=C(F_1(x_1),\cdots,F_p(x_p)).
\]

Most of copulas belong to a class of copulas, specified by the type of dependence it models. For instance, among the upper tail dependence, one could cite the the family of Gumbel copulas~\cite{nelsen}. Thus, copulas provide a simple and natural way to measure the dependence. However, copulas are not the most widely used tool in practice. To measure the dependence, it is usual to refer to the Pearson's correlation coefficient, that measures the linear dependence among variables. It is especially appropriated to elliptical distributions~\cite{fangkt}, but it may be misleading for other types of distribution. The Spearman's rho is a good alternative to the Pearson's coefficient because it could be adapted to any distribution.
Based on the probability of concordance and discordance of random variables, the Spearman's rho is also well-known as a rank correlation, i.e. the linear Pearson correlation coefficient applied on the rank of observations. The main advantage of this measure is that it does not depend on the marginal distributions, but only on the structure of dependence. Furthermore, it is a copula-based measure of association, i.e. when the dependence between two random variables is modelized by a copula $C$, the Spearman's rho, denoted $\rho^{\textrm{S}}$, admits the following expression,

\[
\rho^{\textrm{S}}=12 \int\!\!\!\int_{[0,1]^2} C(u,v) d ud v-3.
\]

Through the Ishigami function, we study the influence of this coefficient on the estimation of our sensitivity measures. We fix $\rho^{\textrm{S}}$, and we model our dependence by two different copulas. The aim is to know if the dependence may be summarized by the  Spearman's rho in the sensitivity analysis in presence of dependent incomes.

\begin{figure}
\centering
\begin{subfigure}[b]{0.45\textwidth}
\centering
\includegraphics[ width = 0.9\textwidth]{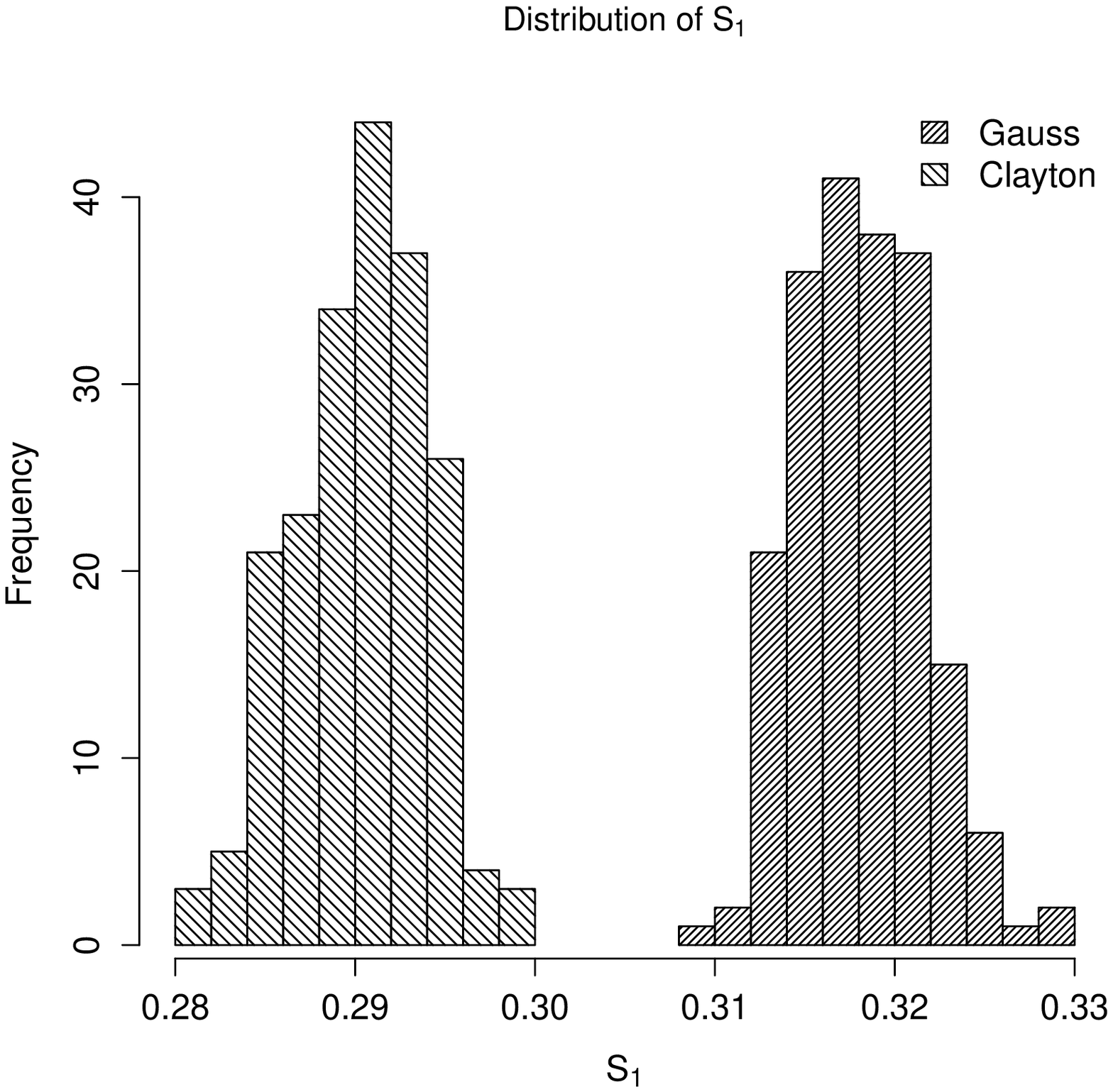}
\captionsetup{format=hang,justification=centering}
\label{figcopulaa}
\end{subfigure} 
\begin{subfigure}[b]{0.45\textwidth}
\centering
\includegraphics[ width = 0.9\textwidth]{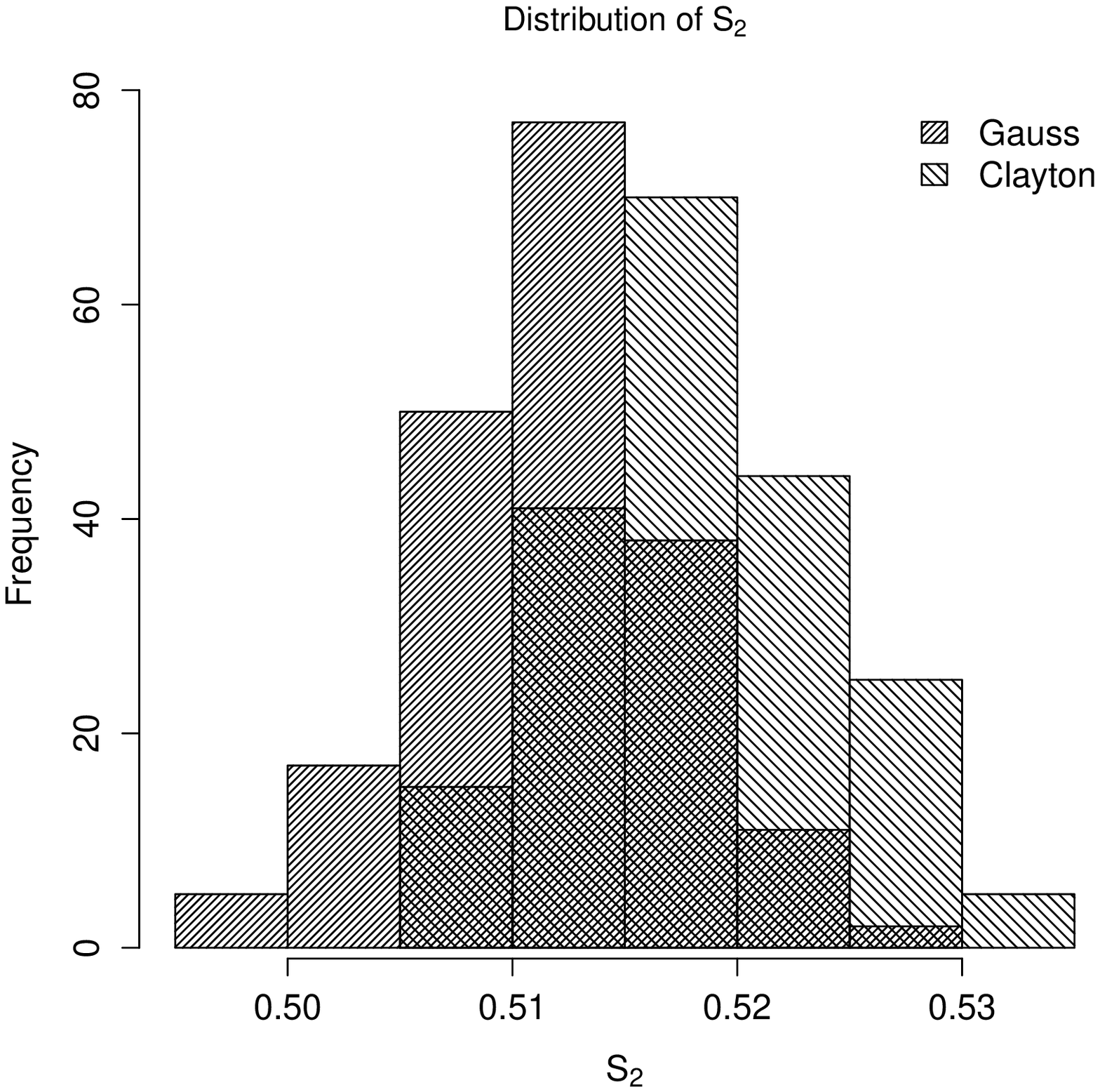}
\captionsetup{format=hang,justification=centering}
\label{figcopulab}
\end{subfigure} 
\\
\begin{subfigure}[b]{0.45\textwidth}
\centering
\includegraphics[ width = 0.9\textwidth]{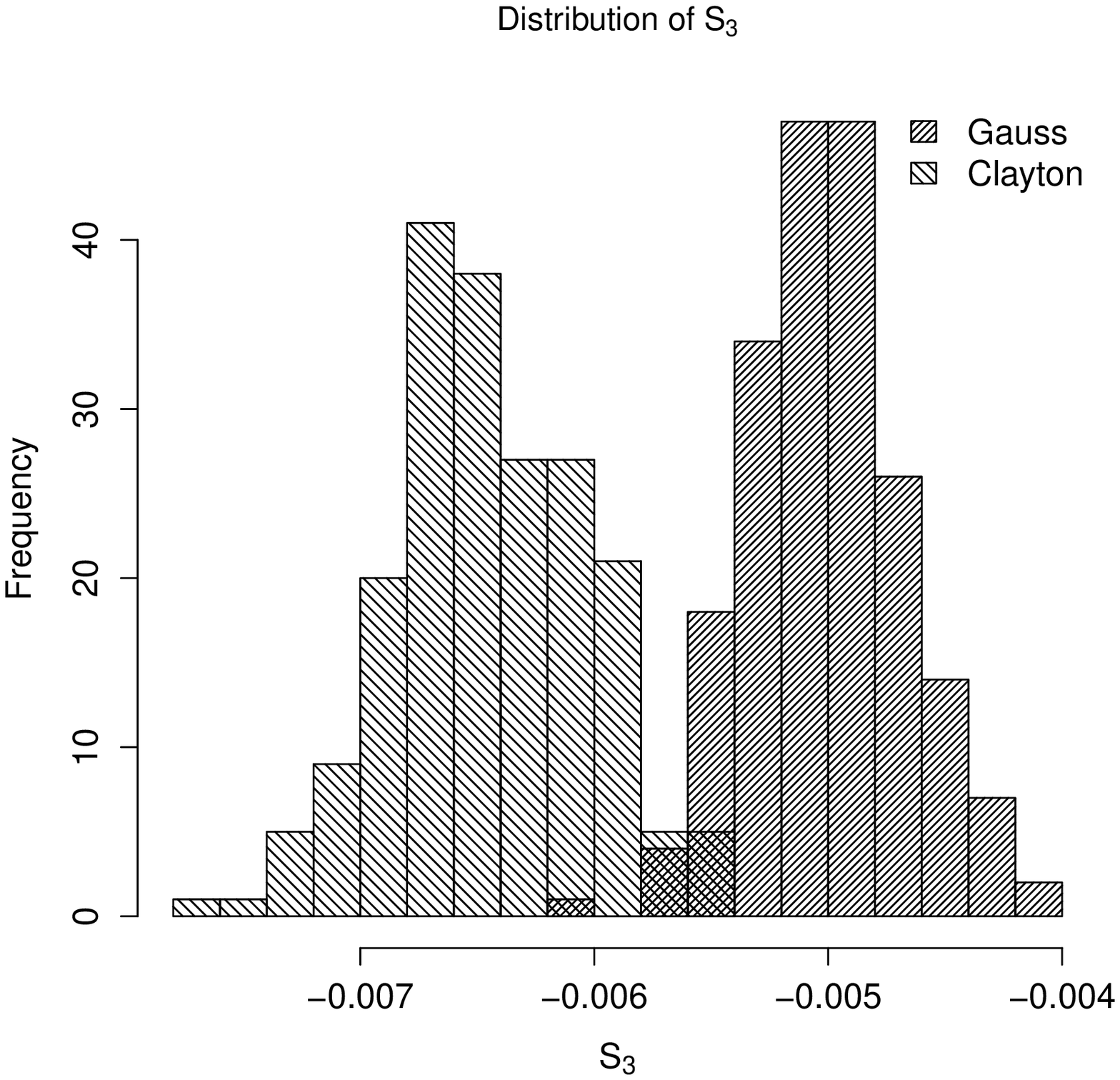}
\captionsetup{format=hang,justification=centering}
\label{figcopulac}
\end{subfigure}
\begin{subfigure}[b]{0.45\textwidth}
\centering
\includegraphics[ width = 0.9\textwidth]{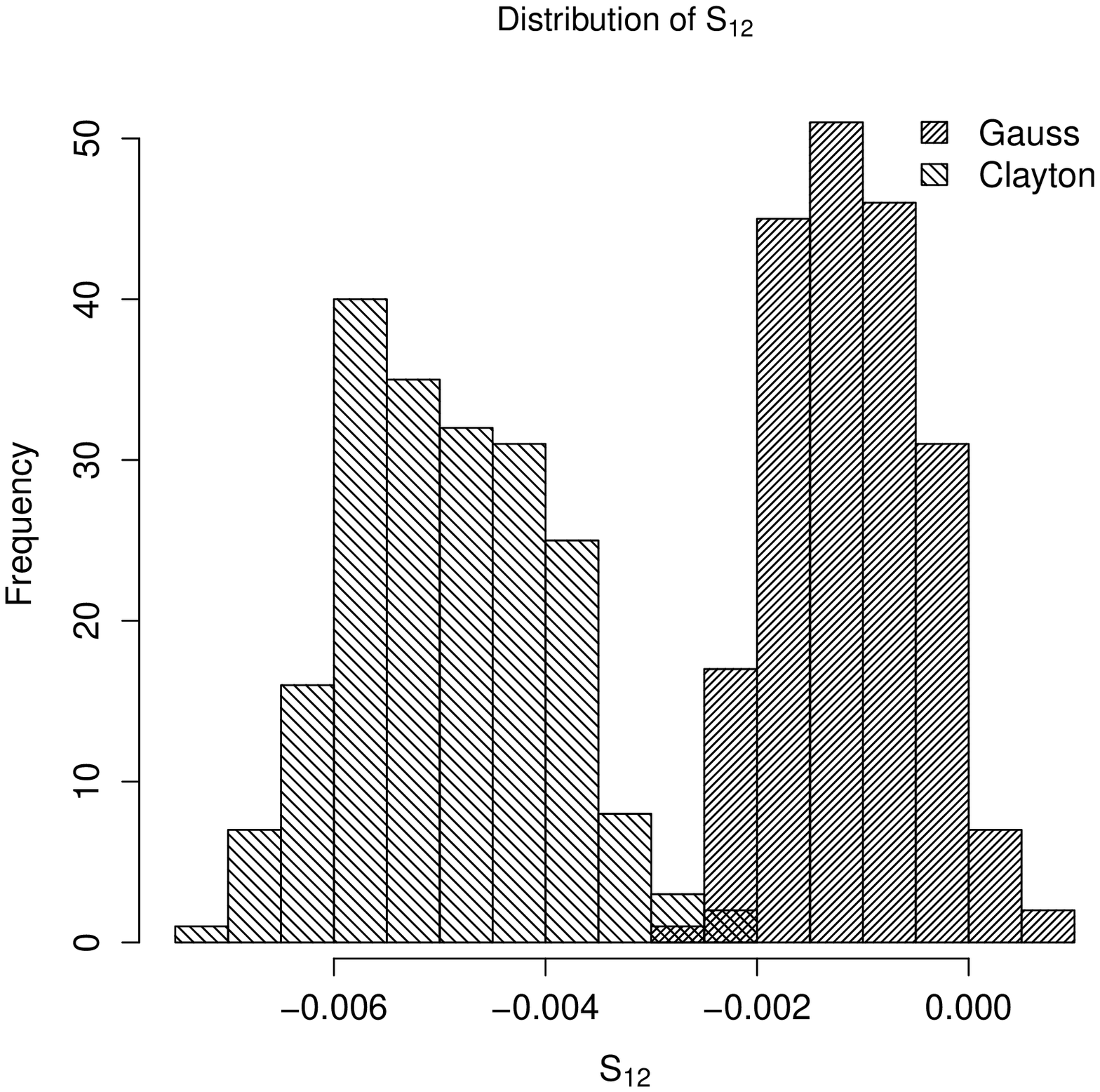}
\captionsetup{format=hang,justification=centering}
\label{figcopulad}
\end{subfigure} 
\\
\begin{subfigure}[b]{0.45\textwidth}
\centering
\includegraphics[ width = 0.9\textwidth]{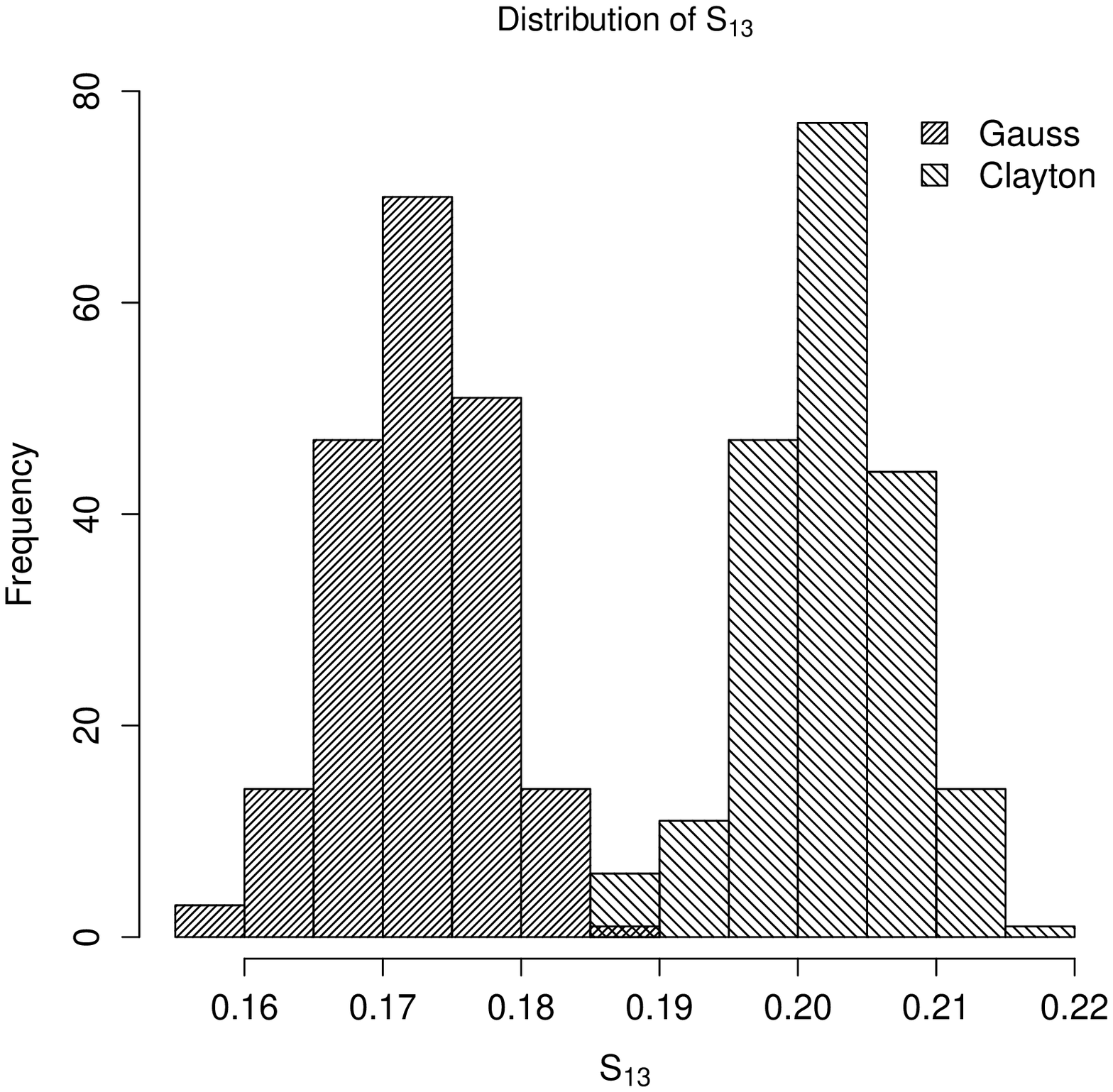}
\captionsetup{format=hang,justification=centering}
\label{figcopulae}
\end{subfigure} 
\begin{subfigure}[b]{0.45\textwidth}
\centering
\includegraphics[ width = 0.9\textwidth]{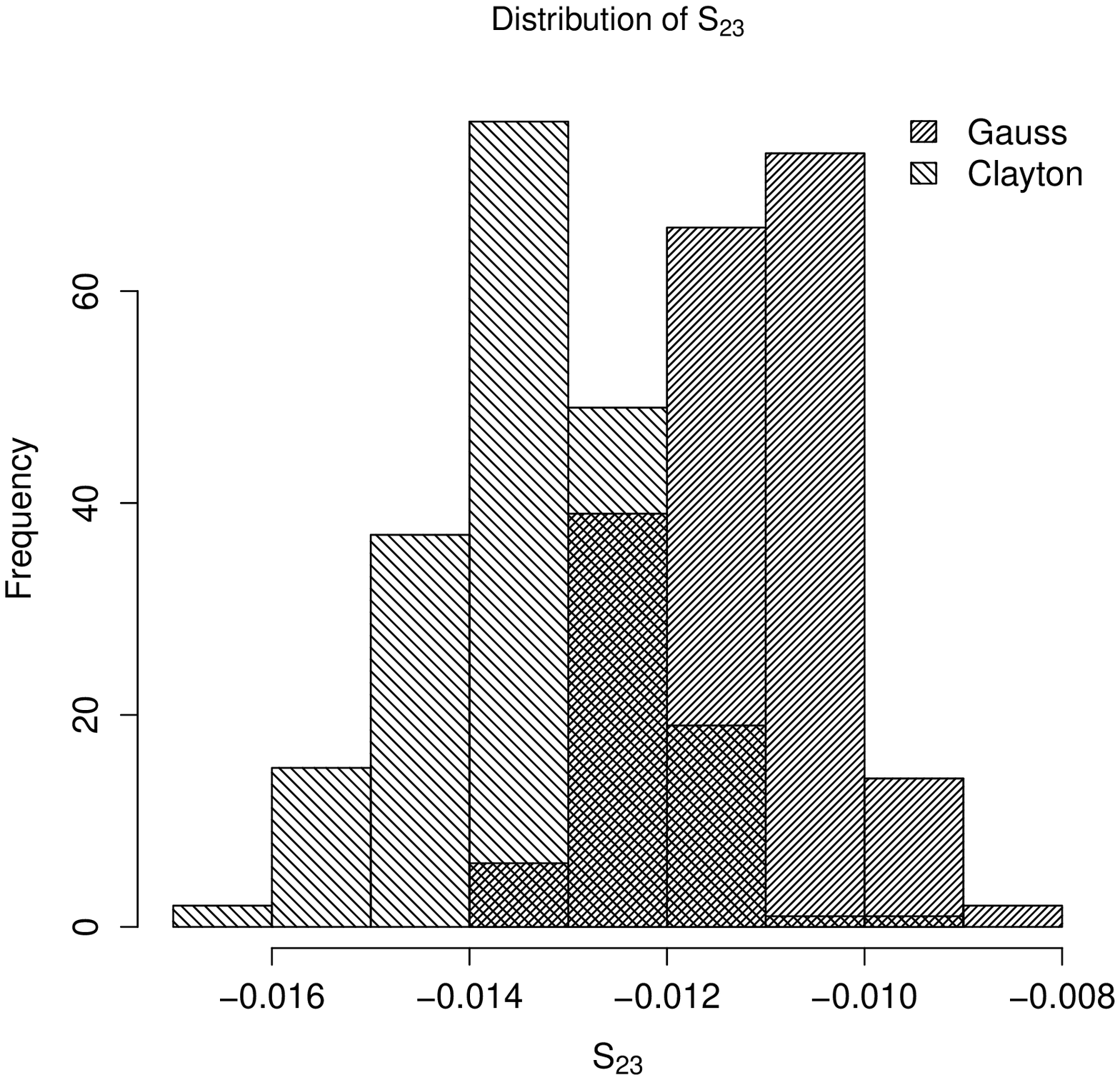}
\captionsetup{format=hang,justification=centering}
\label{figcopulaf}
\end{subfigure}
\caption{Distribution of the sensitivity measures with Gaussian and Clayton copulas}\label{figcopula}
\end{figure}

We assume that each couple of variables $(X_i,X_j)$, $i\neq j$ admits the same Spearman's rho, $\rho^{\textrm{S}}=0.7$. We use the Gaussian copula, and the Clayton copula on uniform marginal distributions over $[-\pi,\pi]$. Further, for a given experimental design set of $n=200$ points, we compare the two dependence structure. Further, the Monte Carlo sample is of size $m=10,000$, and we made $N_s=200$ realizations of $S_{u,m}^f$. 
Figure \ref{figcopula} illustrates the distribution of the indices for both Gaussian and Clayton copulas. \\

Notice that, depending on the type of dependence, we do not obtain the same conclusion. Especially for $S_1$ and $S_{13}$, we notice that the two distributions are disjoint, meaning a significant difference whether we use a Gaussian copula or a Clayton one. 
This shows that it is not enough to only consider a measure of association to model the dependence in sensitivity analysis.

\subsection{Industrial example: the river flood inundation}\label{industrial}

\begin{figure}
\centering
 \includegraphics[width=0.5\textwidth]{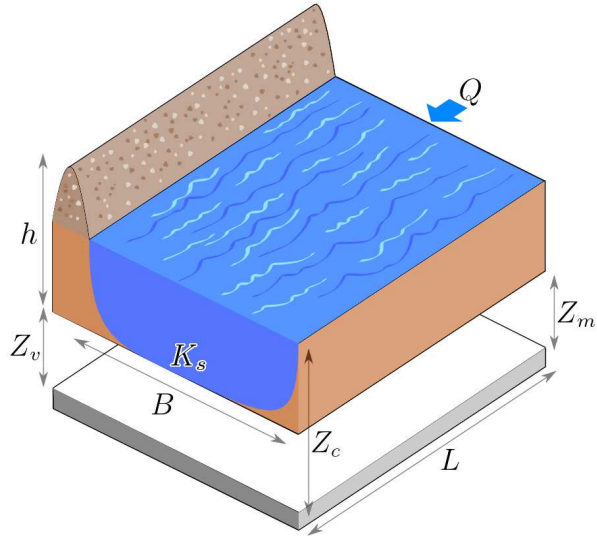}
\caption{The river flood model}\label{floodrep}
\end{figure}

We illustrate our method with the study of a river flood inundation. In this problem, the river flow is compared with the height of a dyke that protects an industrial site~\cite{faivre,rocquigny}. The river flow may lead to inundations that are desirable to avoid. To study this phenomenon, the maximal overflow of the river is modelized by a crude simplification of the 1-D Saint Venant equations, when uniform and constant flow rate is assumed. The model is given by the following expression,

\begin{equation*}
 S=\underbrace{Z_v+h}_{Z_c}-H_d-C_b, \quad h=\left( \frac{Q}{B K_s\sqrt{\frac{Z_m-Z_v}{L}}} \right)^{0.6},
\end{equation*}

where $S$ is the maximal overflow that depends on eight parameters. The river flow and the parameters implied in the model are represented in Figure \ref{floodrep}. These variables are physical and geometrical parameters subject to a spatio-temporal variability or to errors of measurements. Thus, leading a sensitivity analysis in this model has a real interest for this model. The meaning of the incomes and their distribution are given in Table \ref{tableflood}.  

\begin{table}
\centering
\renewcommand{\arraystretch}{1.5}
\begin{tabular}{|c|c|c|}
\hline
Variables & Meaning & Distribution\\
\hline
$h$ &maximal annual  water level & - \\
\hline
$Q$ &maximal annual flow rate & Gumbel $G(1013,558)$ tr. to $[500;3000]$ \\
\hline
$K_s$ &Strickler coefficient & Normal $N(30,8)$ tr. to $[15,+\infty[$ \\
\hline
$Z_v$ &river downstream level & Triangular $T(49,50,51)$ \\
\hline
$Z_m$ & river upstream level & Triangular $T(54,55,56)$   \\
\hline
$H_d$ &dyke height & Uniform $\mathcal U([7,9])$ \\
\hline
$C_b$ & bank level & Triangular $T(55,55.5,56)$\\
\hline
$L$ &length of the river stretch & Triangular $T(4990,5000,5010)$ \\
\hline
$B$ & river width & Triangular $T(295,300,305)$\\
\hline
\end{tabular}
\caption{Description of inputs-output of the river flood model (tr. to=truncated to)}\label{tableflood}
\end{table}

\begin{figure}
\centering
\begin{subfigure}[b]{0.45\textwidth}
\centering
\includegraphics[ width = 1\textwidth]{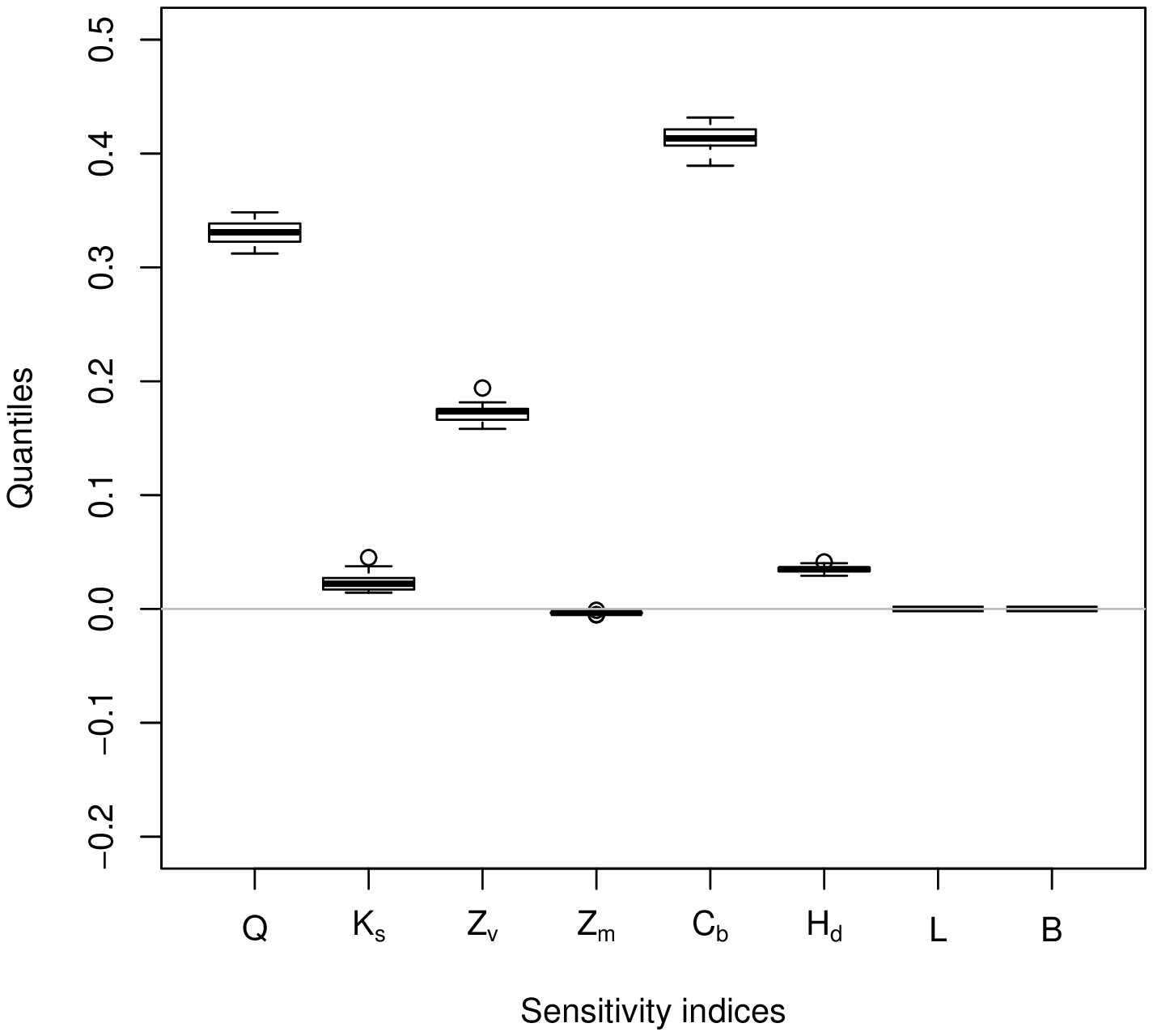}
\captionsetup{format=hang,justification=centering}
\caption[hang,center]{GHOGS procedure}
\label{figcruea}
\end{subfigure}
\begin{subfigure}[b]{0.45\textwidth}
\centering
\includegraphics[ width =1 \textwidth]{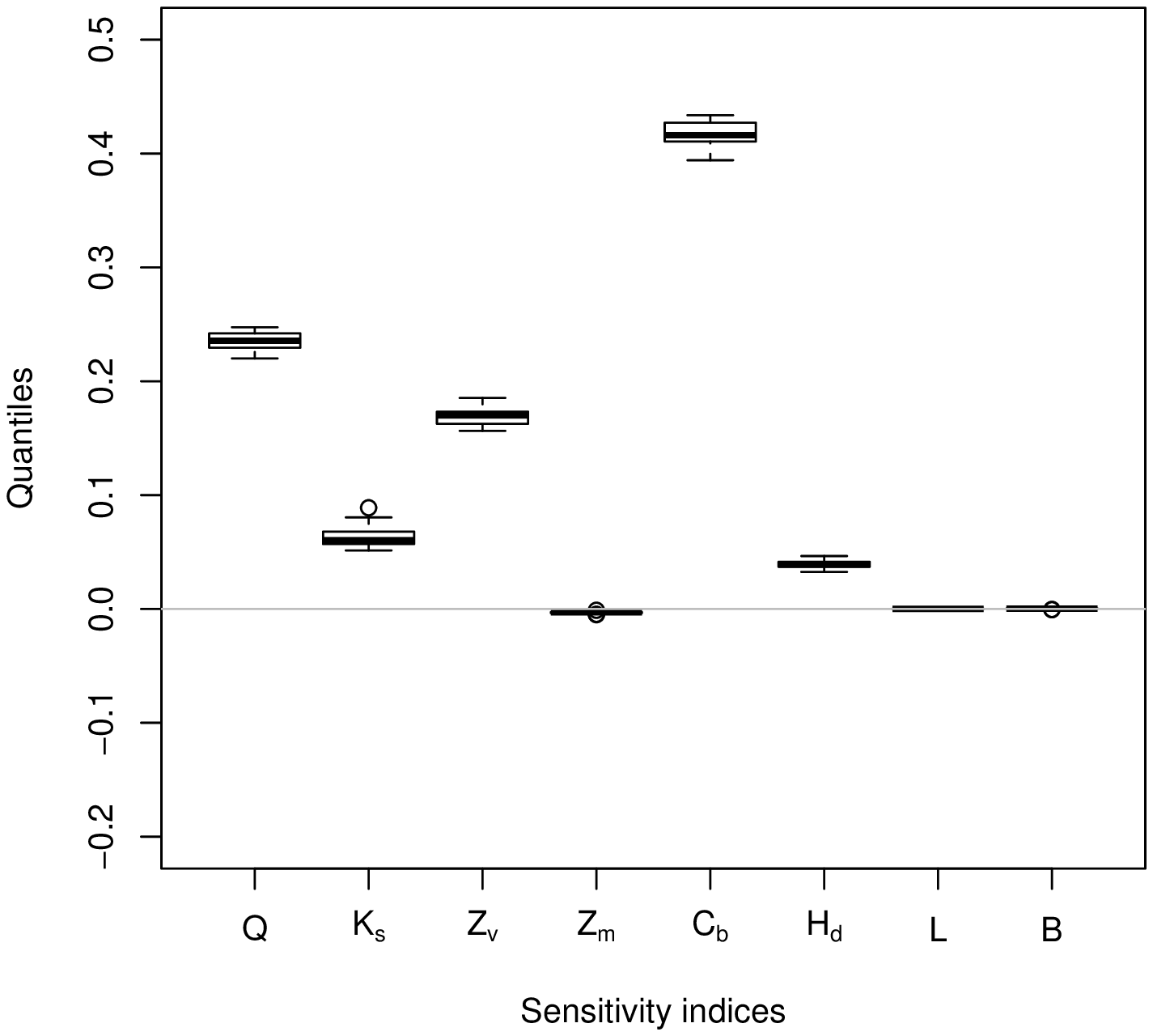}
\captionsetup{format=hang,justification=centering}
\caption[hang,center]{GP modelling}
\label{figcrueb}
\end{subfigure}
\caption{Sensitivity indices estimation with the GHOGS method (\subref{figcruea}) and the GP modelling (\subref{figcrueb})}\label{figcrue}
\end{figure}

In this study, we assume that $(Q,K_s)$ is a correlated pair, with correlation coefficient $\rho=0.5$.
This correlation is admitted in real case, as we consider that the friction coefficient increases with the flow rate.
Also, $(Z_v,Z_m)$ and $(L,B)$ are assumed to be dependent with the same Pearson coefficient $\rho=0.3$, because data are supposed to be simultaneously collected by the same measuring device. As for $C_b$ and $H_d$, they are supposed to be independent. \\
We take a first sample of $n=200$ observations, and a Monte Carlo sample of size $m=5000$. Further, we generate $N_s=100$ realizations of the first order sensitivity indices. We then consider the mode of the probability density estimate of these realizations. We repeat the procedure $100$ times to obtain a Monte Carlo error for each sensitivity index. Further, we compare our result to the generalized sensitivity indices defined in~\cite{chastaing2}, built from a functional decomposition, called \emph{hierarchical} decomposition. The estimation of these last indices is based on a regression approach, and a recursive procedure~\cite{chastaing2}. Further, as this procedure suffers from the curse of dimensionality, a greedy algorithm is adopted to select a sparse number of informative components. This other strategy will be called the GHOGS (for Greedy Hierarchical Orthogonal Gram-Schmidt) strategy. The comparison with our methodology is given by Figure \ref{figcrue}.\\

Through the result, we observe for both decompositions the same phenomena for the last six inputs ($Z_v$, $Z_m$, $H_d$, $C_b$, $L$ and $B$).
The width ($B$) and the length ($L$) of the river are not influent parameters in the model. 
Also, for both decompositions, the dyke height is the most contributive variable in the global variability.
Moreover, the bank level ($C_b$) has a negligible impact on the model  output and its contribution has the same order of magnitude for both analyses.
The main difference between the two procedures is  the contribution of the correlated pair $(Q,K_s)$.   In the GHOGS procedure, we observe that the flow rate $Q$ is highly contributive with respect to the Strickler coefficient $K_s$. This contribution is less important for the Gaussian processes approach while the one of $K_s$ is slightly larger.
Furthermore, we note that the sum of contributions for the pair $(Q,K_s)$ is similar for the two analyses.

\begin{figure}
\centering
\begin{subfigure}[b]{0.45\textwidth}
\centering
\includegraphics[ width = 1\textwidth]{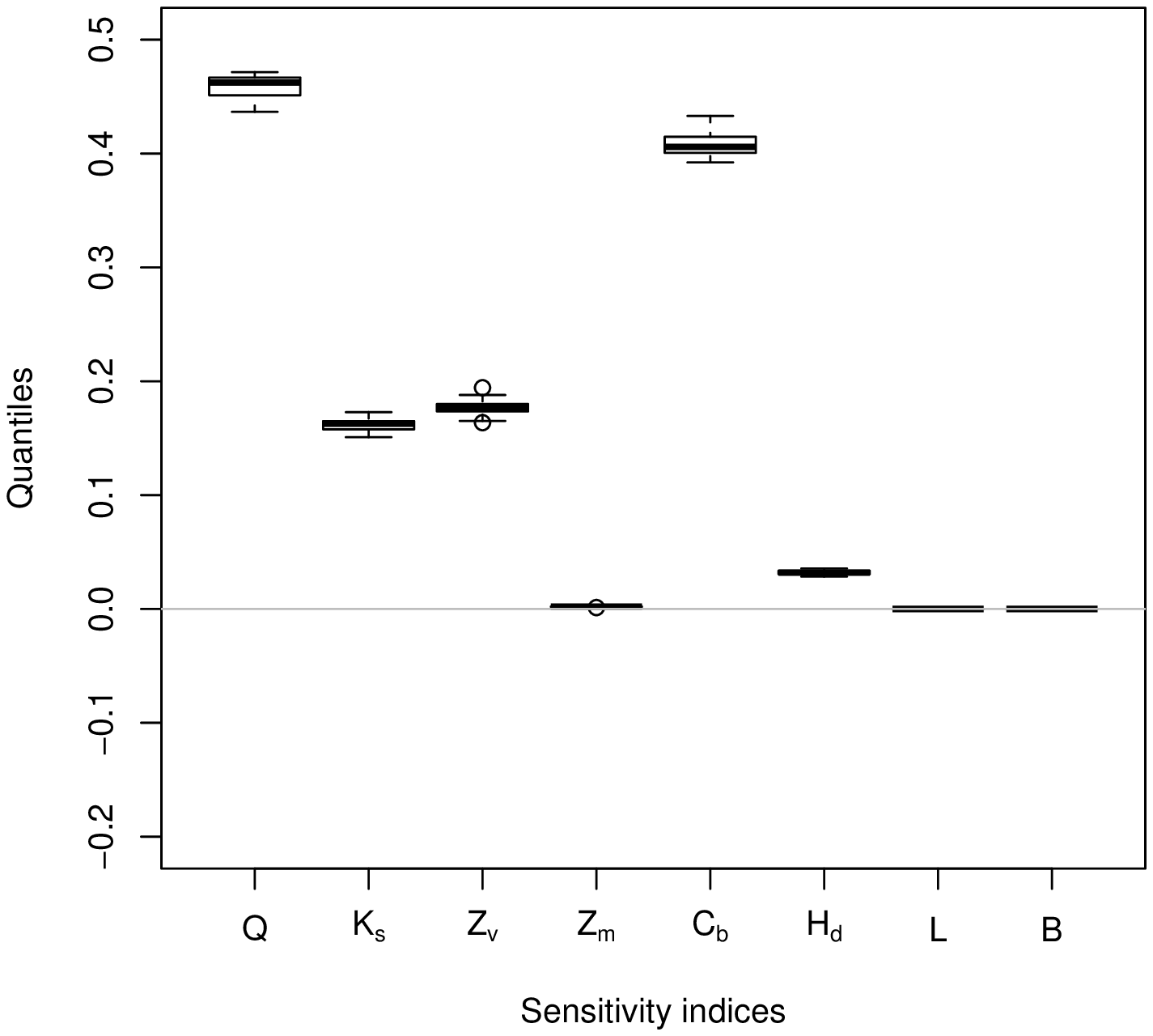}
\captionsetup{format=hang,justification=centering}
\caption[hang,center]{Variance term for the GHOGS procedure}
\label{figcrue2a}
\end{subfigure}
\begin{subfigure}[b]{0.45\textwidth}
\centering
\includegraphics[ width =1 \textwidth]{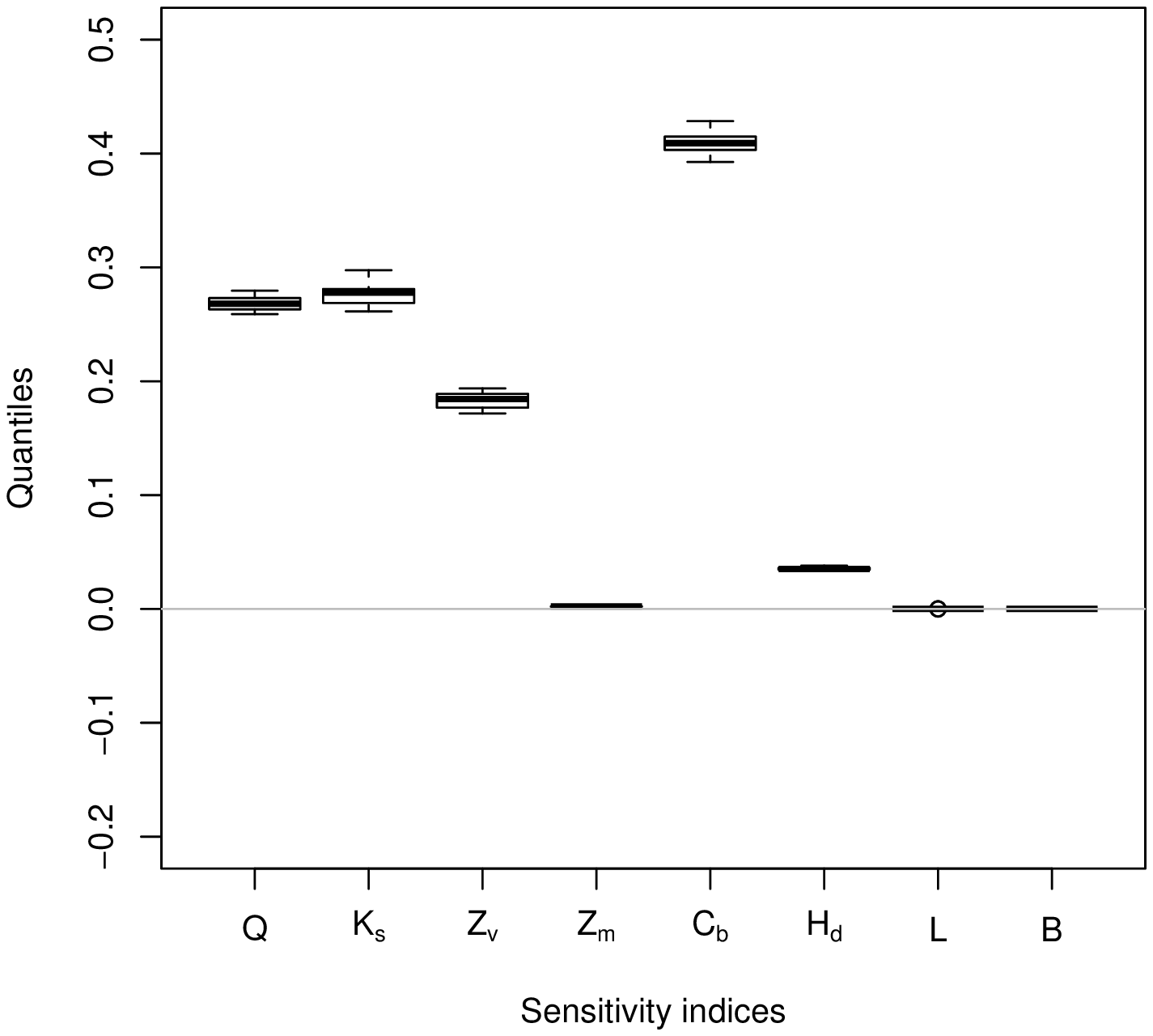}
\captionsetup{format=hang,justification=centering}
\caption[hang,center]{Variance term for the GP modelling}
\label{figcrue2b}
\end{subfigure}
\begin{subfigure}[b]{0.45\textwidth}
\centering
\includegraphics[ width = 1\textwidth]{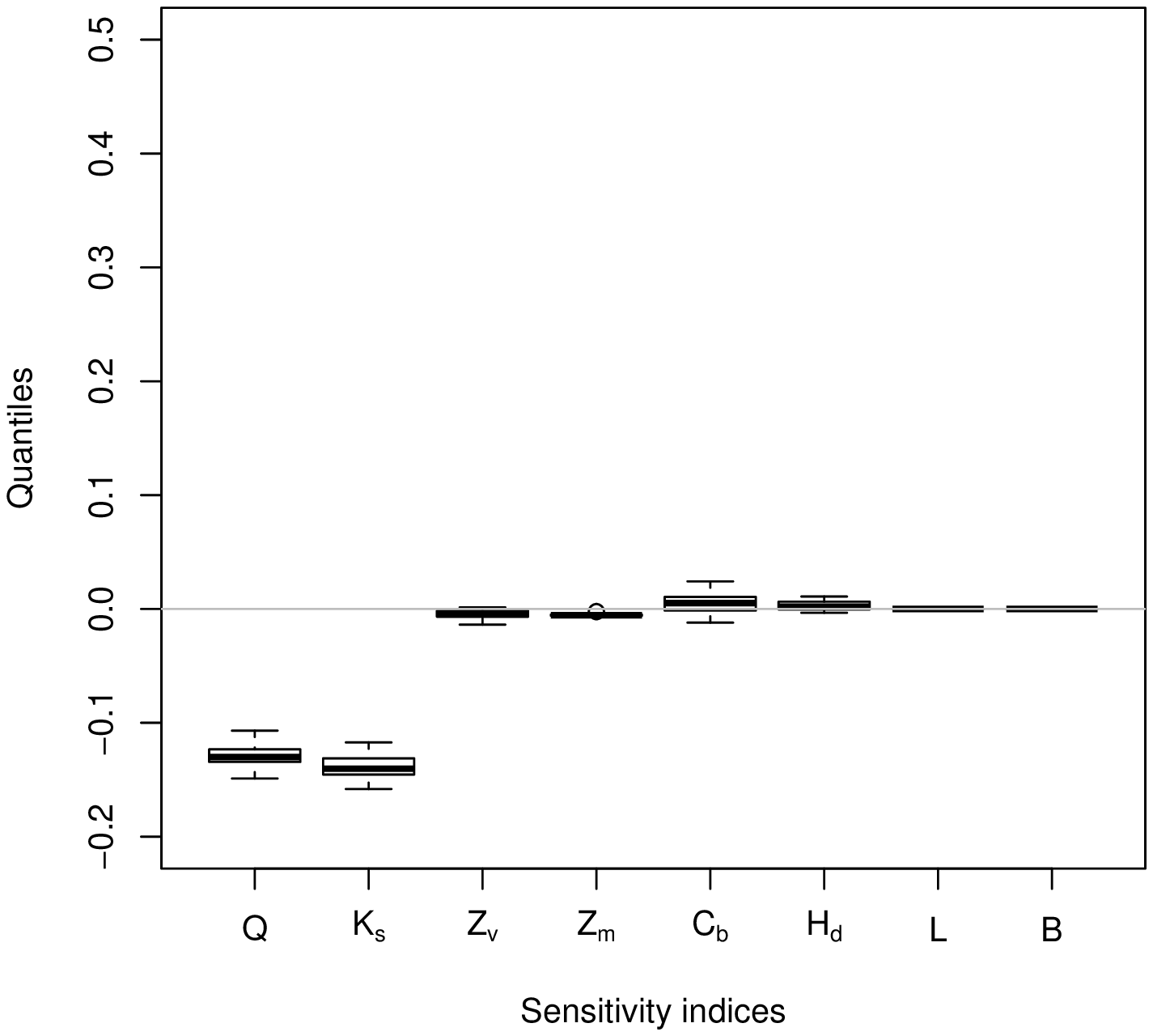}
\captionsetup{format=hang,justification=centering}
\caption[hang,center]{Covariance term for the GHOGS procedure}
\label{figcrue2c}
\end{subfigure}
\begin{subfigure}[b]{0.45\textwidth}
\centering
\includegraphics[ width =1 \textwidth]{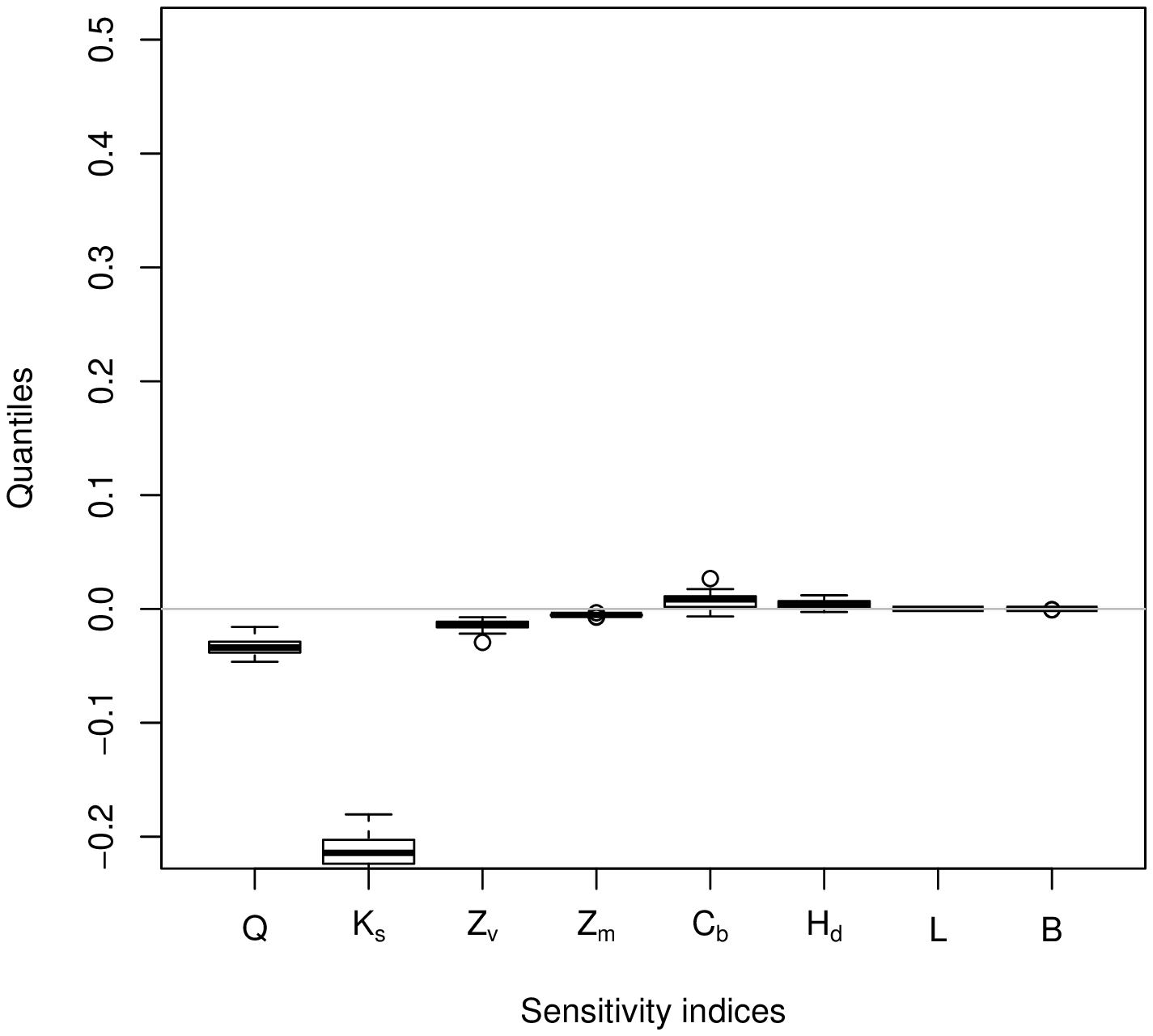}
\captionsetup{format=hang,justification=centering}
\caption[hang,center]{Covariance term for the GP modelling}
\label{figcrue2d}
\end{subfigure}
\caption{Variance (\subref{figcrue2a}) \& (\subref{figcrue2b})  and covariance (\subref{figcrue2c}) \& (\subref{figcrue2d})   terms for the  Sensitivity indices estimation with the GHOGS method (\subref{figcrue2a}) \& (\subref{figcrue2c})  and the GP modelling  (\subref{figcrue2b}) \& (\subref{figcrue2d}) }\label{figcrue2}
\end{figure}

We see in (\ref{indexmetamodel}) that the sensitivity measure is  decomposed into a sum of a variance term $V[ f^n_u(\XX_u)]/V[ f^n(\XX)]$ and a covariance term $\cov[  f^n_u(\XX_u),f^n_{u^c}(\XX )]/V[ f^n(\XX)]$.  The same type of decomposition is  present  in the index provided by the  GHOGS procedure~\cite{chastaing2}. The variance term represents the main contribution of the inputs without the dependence part. The covariance term represents the contribution of the dependence to the index. We represent for the two methods the estimated variance and the covariance parts in Figure \ref{figcrue2}. We observe that the GP and the GHOGS behave differently, as we are not faced to the same decomposition. In the GP approach, the model tends to balance the main contribution, whereas the GHOGS is more discriminant. The covariance contribution is the same in the GHOGS procedure for $(Q,K_s)$, which seems reasonable as it is estimated from a hierarchically orthogonal decomposition~\cite{chastaing}. However, we observe a significant difference between the covariance part of $Q$ and the one of $K_s$, that may be due to the fact that we measure $\cov(f^n_{K_s},f^n_{K_s^c})$. This last term implies a large sum of terms that may weaken the main contribution, and that lead to a negative covariance contribution.


\section{Conclusions}

Through this work, we propose a solution for dealing with complex computer codes in presence of dependent input variables in the model. The definition of a variance-based sensitivity index aims at quantifying the contribution of a (group of) variable(s) in the model and can be decomposed as a sum of ratio of variances, interpreted as the main contribution, and a ratio between covariance terms and the global variance, interpreted as the contribution due to the dependence. The attractive side of such methodology is to be able to quantify the uncertainty of the sensitivity measure, and thus to compute confidence intervals for each estimation. The question about the choice of the ANOVA kernel has not been raised in this work, as this choice may have a strong influence on the values of the sensitivity indices. This remains an open problem.



\newpage

\bibliographystyle{apalike}
\bibliography{biblio_GP}

\end{document}